\documentclass[11pt]{amsart}
\usepackage{latexsym,amsmath,amssymb,epsfig,graphicx}
\usepackage{caption}
\usepackage{cellspace}
\usepackage[utf8]{inputenc}
\usepackage[greek.polutoniko,english]{babel}
\usepackage{wasysym}
\usepackage{wrapfig}
\usepackage{float}

\newtheorem{prop}{Proposition}

\usepackage[left=1.4in,right=1.4in,top=.95in]{geometry}

\begin{document}
\title[The Eclectic Content and Sources of Clavius's \emph{Geometria Practica}]{The Eclectic Content and Sources of Christopher Clavius's \emph{Geometria Practica}}

\author{John B. Little}
\address{Department of Mathematics and Computer Science,
College of the Holy Cross, Worcester, MA 01610}
\email{jlittle@holycross.edu}

\date{\today}

\maketitle
\begin{abstract}
We consider the \emph{Geometria Practica} of Christopher Clavius, S.J., a suprisingly eclectic and comprehensive
textbook of \emph{practical geometry}, whose first edition appeared in 1604.  Our
focus is on four particular sections from Books IV and VI where Clavius has either used his sources in an interesting way or where
he has been uncharacteristically reticent about them.  These include the treatments of 
Heron's Formula, Archimedes' \emph{Measurement of the Circle}, four methods for constructing two mean proportionals 
between two lines, and finally an algorithm for computing $n$th roots of numbers.
\end{abstract}

\section{Introduction}

\subsection{Clavius}

Christopher Clavius, S.J. (1538--1612) was certainly the preeminent Jesuit mathematician of his era and an 
important mathematical astronomer.\footnote{See \cite{Knobloch1988, Baldini1983, Baldini2003}.  
The original documentary sources for many of the facts mentioned here are discussed by Knobloch and Baldini.}  
 He was admitted into the Society of Jesus in Rome in 1555, studied at the University of Coimbra in Portugal for two 
 years, then returned to Rome where he completed his studies.  In comments to some of his younger Jesuit colleagues, 
 he claimed that he was largely self-taught in mathematics.  This is sometimes seen as doubtful since 
 the well-known mathematician Pedro Nunes 
 (1502--1578) was active at Coimbra during Clavius's time there.   However, while there 
 might be traces of Nunes's influence in some of Clavius's more algebraic works, no direct evidence of contact is known. 
 Starting from 1563 through the end of his life (except for a short assignment in Naples in 1595--1596), Clavius 
 served as professor of mathematics at the Jesuit \emph{Collegio Romano}.  At the start of this time, he taught the 
 regular mathematics curriculum and led an ``Academy'' in which exceptionally able and energetic students
could pursue the study of mathematics beyond the basics.  Around the time of his sojourn in Naples, 
he essentially retired from regular teaching and devoted himself primarily to writing and 
mentoring the mathematicians of the Academy.  Many of the talented
Jesuit mathematicians of the generation around 1600 (Christoph Grienberger, Odo van Maelcote, Gr\'egoire de St.~Vincent, Paulo Lembo, Paul Guldin, Orazio Grassi, and others) passed through this Academy and some stayed on at the
\emph{Collegio} as professors of mathematics.   

Clavius was fundamentally a commentator, expositor, and evaluator of the mathematical
work of others, not primarily as an original mathematical researcher in the modern sense.  His mathematical outlook was essentially conservative and grounded firmly in the geometry of the \emph{Elements} of Euclid, with some ``excursions'' into parts of algebra and what we would call discrete mathematics.  Yet his view of the subject was broad enough to acknowledge \emph{both}
the certainty of mathematical knowledge due the subject's reliance on strict standards of proof \emph{and} the utility of mathematics for understanding the physical 
world.\footnote{This is expressed most explicitly in 
Clavius's essay \emph{In disciplinas mathematicas prolegomena} (\emph{Prolegomena} on the mathematical
disciplines) included in Volume I of the \emph{Opera Mathematica}, \cite{Clavius1611}.  Clavius sees mathematics as \emph{intermediate}
between metaphysics and natural philosophy, an idea that traces back at least to Proclus's Commentary on 
Book I of Euclid's \emph{Elements}, \cite{Proclus}, a text Clavius mentions several times. See, for instance,
\cite[Chapter~1]{Rommevaux2005}.}
Besides his teaching and mentoring of younger mathematicians, a large portion of Clavius's energies were 
devoted to the production of numerous influential textbooks or source books for the teaching of a wide range of mathematical subjects.  
These included his extensively augmented edition of the \emph{Elements} of Euclid (first edition 1574), 
the \emph{Epitome arithmeticae practicae} (Summary of Practical Arithmetic, 1583), this \emph{Geometria Practica} 
(Practical Geometry, first edition 1604), and the 
\emph{Algebra} (1608).  Clavius also wrote a well-known commentary on the \emph{Sphere} of Sacrobosco, books on the astrolabe and the construction of sundials, and more elementary treatments of plane and spherical triangles. 
 His collected mathematical works (the \emph{Opera Mathematica}) were published in five volumes starting in 1611.\footnote{These have been digitized, see: {\tt https://clavius.library.nd.edu/mathematics/clavius}.}   

\subsection{The Geometria Practica}

The 1604 edition of the \emph{Geometria Practica} was printed by the shop of Luigi Zanetti in 
Rome; just two years later, a second edition was produced by the printshop of Johann Albin in 
Mainz in 1606.   The 1604 edition is slightly longer because of a different page format.  
However, there are no substantial differences between the texts.  Moreover, very similar (but not identical) 
woodcut figures were used in both editions, so the overall 
appearance does not differ significantly.    The version of the \emph{Geometria Practica} included in the \emph{Opera 
Mathematica} contains some corrections of typographical and mathematical errors in the previous editions, 
an expanded discussion of the quadrant constructed in Book I, and some other relatively minor additions. 
In this essay, page numbers refer to the page in the 1606 edition.  

While a number of scholars have written on aspects of the book, as a whole, Clavius's \emph{Geometria Practica} 
has apparently received less scholarly attention than his edition of Euclid's
\emph{Elements}.\footnote{This is perhaps a reflection of certain derogatory attitudes toward ``applied'' or ``practical''
mathematics in general.}     Yet this work has unexpected and surprising features.   Clavius's 
account certainly stands out in several ways within the whole genre of \emph{practical geometry}
texts,\footnote{See \cite{Raynaud2015} and the Introduction to \cite{Hugh} by the translator, F. Homann, S.J.}
and these features form the major reasons this work remains interesting from the historical 
perspective.  They have also  furnished the main motivations for the author of this essay to undertake a translation of the 
entire \emph{Geometria Practica} from the original Latin into English using the 1606 second edition.\footnote{This translation is available at {\tt CrossWorks}, the online 
faculty and student scholarship repository maintained by the Library of the College of the Holy Cross. 
All quotations of passages from the \emph{Geometria Practica} in English are taken from this translation.  The original Latin text from the 
1606 edition of the \emph{Geometria Practica} will be provided in footnotes, for purposes of comparison.}

First, the \emph{eclecticism}--the sheer range of different types of topics that fall under the 
category of \emph{practical geometry} or allied areas for Clavius and that make it into this book--is remarkable.  
In his Preface, Clavius discusses previous works in the \emph{practical geometry} genre, and how
he wants his work to stand out from the others:
\begin{quote}
... [M]any erudite men
have pursued all of its [\emph{i.e.}~practical geometry's] parts with accurate and diligent writing.   Among them, 
Leonardo Pisano [``Fibonacci''], Brother Luca Paccioli, Nicolo Tartaglia, 
Oronce Fin\'e,  Girolamo Cardano and others have demonstrated preeminence and 
flourished to exceptional praise.  But I would judge Giovanni Antonio Magini one of the first
in mathematical excellence.  He has taught so much about the measurement of lines, 
and treated this subject so fully, so systematically and with such perspicacity, that
he seems to have snatched away, not only the standing from those who wrote before him, but also the hope of equal, 
let alone greater, glory from those coming after him.  But truly,  Magini concerned himself only with this one 
part of this subject, and the others, although they undertook to present all of those parts, have left out much 
in writing their books.  I decided, if possible, to complete the subject, so that \emph{whatever has been profitably 
handed down by others or found by myself in practical geometry is enclosed within the circle of one work}.\footnote{
...  \& multos, \& eruditos viros ... , qui partes illius omnes accurata, \& diligenti scriptione persecuti sunt:
Inter quos, vt Leonhardus Pisanus, Frater Lucas Pacciolus, Nicolaus Tartalea,
Orontius, Cardanus, aliique praecipuas obtinuerunt: ita eximia in caeteris 
laude floruerunt.  Primas tamen adiudicarim Io. Antonio Magino praestanti 
Mathematico; qui tametsi tantum linearum dimensiones docuit, ea tamen copia,
doctrina, perspicacitate cuncta tradidit, vt locum non modo iis, qui ante 
scripserunt, sed spem posteris aequalis gloriae, ne dum maioris, ademisse
videatur.  Verum quoniam \& hic de vnica tantum parte fuit sollicitus: \& alii, quamuis aggressi
omnia, multa tamen inter scribendum praeterierunt: decreui, si qua possem, 
perficere: vt, quicquid vtiliter in Geometria practica ab aliis traditum, \`a me 
etiam inuentum est, vnius operis gyro clauderetur. } (\cite[Preface]{Clavius1606}; emphasis added) 
\end{quote}

Thus, in this work, Clavius aims for \emph{completeness} within the subject of practical geometry
as understood by his contemporaries.  
The following rough outline of topics will demonstrate how wide-ranging and encyclopedic this book truly is:\footnote{At a 
higher degree of granularity, the complete list of chapter headings and propositions that serves as the 
table of contents is even more evidence here.} 
\begin{itemize}
\item[-] Book I:  Construction of a proportional compass and quadrant for measuring lengths 
and angles; summary of elementary plane trigonometry
\item[-] Book II:  Measuring lengths, heights, and depths with the astronomical quadrant\footnote{Discussions of 
problems similar to those considered here can be seen in almost all practical geometry books.   
As Raynaud  points out, \cite[p. 15]{Raynaud2015}, these are part of a long and surprisingly stable tradition 
with connections to propositions 19--22 from Euclid's \emph{Optics}.}
\item[-] Book III:  A parallel discussion of measuring lengths, heights, and depths with the 
geometer's square; other methods for the same sorts of problems
\item[-] Book IV: Measuring areas of plane regions, including an augmented translation of
of Archimedes' {\it Measurement of the Circle}, and quoting from other works of Archimedes including
the \emph{Quadrature of the Parabola}
\item[-] Book V: Measuring volumes of solid bodies, with extensive quotations of 
results from Archimedes' works {\it On the Sphere and Cylinder}, and \emph{On Conoids and Spheroids}.
\item[-]  Book VI:  Geodesy, that is, the division of rectilinear surfaces of whatever
sort, either by lines drawn through some point, or by parallel lines;\footnote{Clavius's 
sources and his treatment of this topic are discussed
by E. Knobloch in \cite{Knobloch2015}.}   how 
plane or solid figures are increased or decreased in a given ratio;  several methods for 
finding two mean proportionals between two given lines selected from the commentary on Archimedes' 
{\it On the Sphere and Cylinder} by Eutocius of Ascalon and  Pappus's \emph{Mathematical
Collection};  finally, an algorithm for extracting all sorts of roots by hand calculations 
\item[-]  Book VII:  Isoperimetric figures and questions (drawing on material in Pappus and the commentary on 
Ptolemy's \emph{Almagest} by Theon of Alexandria), together with an appendix on 
the problem of squaring the circle via the \emph{quadratrix} curve of Hippias, drawing from Pappus 
\item[-]  Book VIII: Various geometric theorems and constructions that Clavius says can be 
used to build mathematical power in problem-solving--several of these are drawn from 
Pappus's \emph{Mathematical Collection}, including one discussing the trisection of 
general angles using the \emph{conchoid} curve of Nicomedes.  A table of squares and 
cubes of all $N \le 1000$ is included at the end, together with a discussion of how the table can be
extended using facts about the first and second differences of the sequences of squares and cubes,
with applications to extraction of square and cube roots.
\end{itemize}

As is true in all of his other works, Clavius also has clarity of exposition as a second main goal.  
A striking example of this commitment to completeness and clarity of exposition is 
Clavius's treatment of Archimedes' \emph{Measurement of the Circle} in Chapter 6 of Book IV, which will 
be examined in detail in \S \ref{ArchimedesMeasurementCircle}.

The second feature that has seemed surprising to this author is the resolutely \emph{dual theoretical and practical
focus} of much of this text on \emph{practical geometry}.\footnote{As Knobloch writes,
``... son approche d\'emontre les limites d'une division trop tranch\'ee entre g\'eom\'etrie
pratique et g\'eom\'etrie savante," \cite[p. 60]{Knobloch2015}.  That is, 
``... [Clavius's] approach shows the limitations of a too-definite distinction between
practical geometry and theoretical geometry.''  This applies to almost every section of the 
\emph{Geometria Practica}, not just the discussion of geodesy.}  
The practical side is signaled immediately
in the author's Preface where Clavius discusses his motivation for writing the book.  After saying
that his experience as a teacher has taught him that most students work and learn best when they understand 
that what they are learning will prove to be useful,\footnote{Et ver\`o cum perpetua multorum annorum 
experientia compererim, admodum paucos esse, qui non in Mathematicis 
exerceantur eo consilio, vt quae didicerint, ad aliquem vsum trahant.}
Clavius addresses how the contents of this book may 
find uses in the real world:
\begin{quote} 
For of course as long as the methods
by which we must make measurements to understand the lengths of fields, the heights of mountains, the 
depths of valleys, and the distances between all locations are presented,
it is clear to anyone (in my opinion) how much that is of use in the construction of buildings, in 
agriculture, in the design of weapons, in the contemplation of the stars, and in all other 
arts and disciplines, can flow from the study of these things.\footnote{Etenim 
dum certa ratio traditur, qua camporum longitudines, altitudines montium, 
vallium depressiones, locorum omnium inaequalitates inter se, \& interualla 
deprehendere metiendo debeamus: cuilibet liquet, vt arbitror, quantum commodi,
vtilitatisque substructioni aedificiorum, cultui agrorum, armorum tractationi, 
contemplationi siderum, aliisque artibus, \& disciplinis ex horum cogitatione 
manare possit.} (\cite[Preface]{Clavius1606})
\end{quote}
Clavius consistently uses numerical examples in many sections and he presents a number of purely calculational 
methods (e.g. the methods for extraction of roots in Book VI and the material on differences of 
squares and cubes at the end of Book VIII).  He discusses the use of different mechanical tools for measurements
and is even willing to countenance
``\emph{mechanical},''  hence necessarily approximate, methods of measurement in geometric diagrams:  
\begin{quote}
No one should be troubled that we have said lines are sometimes to be measured mechanically with an iron chain, or by means of the \emph{instrumentum partium}.\footnote{That is, the proportional compass introduced in Book I.}
 For in this business, especially for fields and farms, this mechanical way of making measurements is wholly admissible, partly since this is the custom among all surveyors, partly because the geometric way is not always possible, but mostly because for the dimensions of farms or other areas it is sufficient to come close enough to the truth that no notable error is made. If someone does not approve of this way of measuring lines, it assuredly and necessarily takes away every possibility of measuring farms or other areas. For in what way is it certain that a given field or figure has known sides, if these have not been explored by some material measurement? If therefore mechanical measurement of lines (not straying far from the truth, as it were) is used by everyone, I do not see why we should think to reject it in measuring lines in figures.\footnote{Neminem autem moueat, aut perturbet, quod rectas dixerimus metiendas esse
nonnunquam mechanice per catenulam aliquam ferream, aut per instrumentum
partium.  Nam in hoc dimitiendi negotio, praesertim in campis, \& agris admittenda
omnino est huiusmodi mechanica linearum dimensio, tum quia apud omnes                  
agrimensores hic mos est: tum quia non semper via Geometrica id praestare                          
potest ; tum vero maxim\`e, quia in dimensionibus agrorum, siue figurarum satis                                
est rem prope verum attingere, dum modo notabilis error non commitatur.  Quod 
si haec dimensio quarundem linearum alicui non probetur, is profecto \`e medio
tollat, necesse est, omnem agrorum, figurarumue dimensionem.  Vnde enim
constat, agrum propositum, vel figuram habere latera cognita, nisi haec ipsa
per mensuram aliquam materialem sint explorata?  Si igitur laterum dimensio
mechanica, tanquam \`a vero parum aberrans, ab omnibus vsurpatur, cur eam in
lineis intra figuras metiendis reijciendam censeamus, non video.}   
(\cite[p.~169]{Clavius1606}])
\end{quote}

There are indeed possible practical applications of many of the more theoretical topics treated in Books 
IV through VIII as well.  Clavius includes a section on methods
used by surveyors in Book IV and a section on measuring volumes of barrels or 
casks at the end of Book V.  However, once he gets past the very basic material in Books I, II, and III,
Clavius's focus seems to shift to developing the mathematical theory \emph{along with} a few practical 
applications, and he usually provides full proofs for the most important results.  

Another hallmark of Clavius's approach and theoretical orientation even within the practical discussions 
is his scrupulous attention to providing reasons for 
almost everything he writes and sources for the material he does not prove in detail.   
This applies even within Books I, II, and III.
Throughout the text, an elaborate system of marginal notes identifies justifications for assertions and 
for the individual steps in proofs or computations.  Over the course of the whole book, the justifications
for the steps in those proofs span almost all of the 13 books of the canonical version of Euclid's
\emph{Elements}, plus the 14th and 15th books added by later authors and included in Clavius's edition
of the \emph{Elements}, as well as some of 
Clavius's other texts, several works of Archimedes,\footnote{Primarily \emph{Measurement of the Circle}, 
\emph{On the Sphere and Cylinder}, and \emph{On Conoids and Spheroids}, but less commonly also
the \emph{Quadrature of the Parabola}.}  and Apollonius's \emph{Conics} (once).

A third feature that is clearly visible in the above outline, but that might be surprising,
is the extent to which Clavius draws on 
Archimedes, Pappus, Claudius Ptolemy, Euclid's \emph{Elements} and Proclus's commentary on Book I,
and works of other ancient and medieval mathematicians and his contemporaries.\footnote{Interestingly enough,
Clavius only refers to the {\it Conics} of Apollonius a handful of times.}   It is worthwhile
to note here that some of the work of the ancient Greeks was just coming back into the European mathematical
mainstream at precisely this time due to the work of humanist scholars such as Federico Commandino and others.
Commandino's Latin translation of the surviving portions of Pappus's \emph{Mathematical 
Collection}, for instance, only appeared in print in 1588.   

\subsection{This essay}

Our plan in this essay is to flesh out this general description of the eclectic content and sources of Clavius's 
\emph{Geometria Practica} by focusing on four particular sections dealing with topics of particular interest.
We have restricted ourselves to parts of the text not covered in detail by other authors.  So for instance,
we have not included a discussion of the Appendix to Book VII giving Clavius's approach to the problem of squaring the circle via the 
\emph{quadratrix} curve\footnote{A more extensive version of this also appears in Clavius's edition of Euclid.} because that is analyzed deeply 
by Bos in \cite[Chapter 9]{Bos}.
Similarly we have not considered the discussion of \emph{geodesy} at the start of Book VI, since Clavius's
approach has been discussed by Knobloch in \cite{Knobloch2015}.   The sections we do discuss are ones
where (in our judgment) Clavius has either used his sources in an interesting way, or he has been
uncharacteristically reticent about those sources.\footnote{As a rule, Clavius is very careful to identify sources, and
it stands out when he does not do so.  
Over the course of this book, the list of authors cited is quite extensive, including (but possibly not limited to)
Apollonius, 
Archimedes, 
Archytas,
Giovanni Battista Benedetti, 
Campanus de Novare, 
Girolamo Cardano, 
Federico Commandino,
John Dee, 
Dinostratus,
Diocles,
Albrecht D\"urer, 
Eratosthenes,
Euclid,
Eutocius of Ascalon,
Oronce Fin\'e, 
Fra\c{c}ois de Foix, Comte de Candale,
Niccolo Fontana (``Tartaglia''), 
Gemma Frisius, 
Marino Ghetaldi,
Christoph Grienberger, 
Hippocrates,
Hypsicles,
Ioannes Pediasimos,
Leonardo Pisano (``Fibonacci''),
Ludolph van Ceulen, 
Mohammad of Baghdad,
Odo van Maelcote,
Giovanni Antonio Magini,
Francesco Maurolico,
Menaechmus, 
Nicholas of Cusa, 
Nicomedes, 
Latino Orsini,
Luca Pacioli,
Pappus, 
Georg Peuerbach,
Proclus, 
Ptolemy,
Joseph Justus Scaliger, 
Sporus, 
Simon Stevin, 
Theon of Alexandra, 
Juan Bautista Villalpando, 
Johannes Werner.   A fuller listing of all the authors cited by Clavius across his whole written output 
is given in \cite{Knobloch1995}.}
  We will look first at the beginning of his discussion of computing 
areas of triangles in Book IV, where Clavius presents what we now call \emph{Heron's formula} before the usual 
method based on finding an altitude of the triangle.  He does not say
anything about his sources there, but by comparing what Clavius says with the development in Leonardo Pisano's
\emph{De Practica Geometrie} some insight may be gained.  Second, we will look at Clavius's treatment of 
Archimedes' {\it Measurement of the Circle} in Book IV.  Here we will see that Clavius has presented essentially a 
complete reworking of the extant Archimedean text incorporating many additional explanatory comments and details 
not found in other versions.  Third, we will consider Clavius's discussion of some of the Greek constructions for finding
two mean proportionals between two given lines.  This involved a very deliberate selection of only a few of the 
methods discussed in the commentary on Archimedes' \emph{On the Sphere and Cylinder} by Eutocius of 
Ascalon and by Pappus in the \emph{Mathematical Collection}.  Finally, 
we consider the discussion of an algorithm for extraction of $n$th roots discussed at the end of Book VI.  
Clavius does not explicitly identify his source here.  But by considering what books would have been available to him, 
and comparing his treatment of extraction of roots with what appears in one of those books, 
 we are able to propose what we believe is a very likely candidate.  This may have been noted before, 
 but if so, we are not aware of it.

\section{Clavius's treatment of ``Heron's formula'' for triangles in Book IV.}

In Chapter 2 of Book IV, Clavius discusses methods for finding the area of a plane triangle.  
He says there are two ways of doing this and he will first present the most accurate or precise one.
He states this as a \emph{rule} or procedure for doing the computation:
\begin{quote}
Let all the sides be added
together in one sum; let each of the sides be subtracted from half of this sum, 
so that three differences between the semiperimeter and the sides are
obtained; finally, let these three differences and the semiperimeter be multiplied
together.  The square root of the number produced will be the area of the 
triangle which is sought.\footnote{Colligantur omnia latera in unam summam:
Ex huius summa semisse subtrahantur singula latera, vt habeantur tres differentiae inter
illam semissem, \& latera singula:  Postremo tres hae differentiae, \& dicta semissis inter
se mutuo multiplicentur.  Producti enim numeri radix quadrata erit area trianguli quaesita.}
(\cite[p.~158]{Clavius1606})
\end{quote}
In modern algebraic terms, the procedure can be collapsed into the single formula 
$$A = \sqrt{s(s - a)(s - b)(s - c)},$$
where $a,b,c$ are the side lengths and $\displaystyle s = (a + b + c)/2$ is the semiperimeter of the triangle.  
This usually goes by the name \emph{Heron's Formula} today, and indeed this is stated
and proved in Proposition I.8 of the \emph{Metrica} of Heron of Alexandria (ca.~10--ca.~70 CE(?)).\footnote{The
Islamic mathematician al-B\={\i}r\=un\={\i} (973--1048) thought that the result was originally proved by Archimedes, and
C.~M.~Taisbak has recently provided a conjectural reconstruction for how Archimedes might have
stated the result.  See \cite{Taisbak2014}.  The statement we have is quite unusual 
for Greek mathematics from the time of Archimedes or earlier because it is not clear what 
geometric significance should be attached to the product of four lengths or the product of 
two areas (which arises in the proof).  In our formula $A = \sqrt{s(s-a)(s-b)(s-c)}$, the quantities
on the right would be interpreted as numbers; in fact Clavius says exactly this at one point in his version of the 
argument. But that is not what Greeks working in the strict Euclidean tradition 
would have done.  Taisbak thinks the Archimedean form
of the statement could have been that the triangle is the mean proportional between two rectangles.  This
is certainly possible. But it must be treated as a conjecture since no known Archimedean text deals with 
questions related to Heron's formula.}
Clavius provides three numerical examples of triangles with integer side lengths, the last of which 
leads to a product $s(s -a)(s - b)(s - c)$ that is not a square.  He then gives a complete, detailed
proof that this does in fact produce the area of the triangle.

Unusually for him, Clavius does not provide an attribution for this result, so several natural questions arise.  First,
what if any source(s) was he drawing on here and why did he not mention it (or them)?  And even before that
question: What source or sources for this result would have been available?
This is an interesting question because in Clavius's time Heron's \emph{Metrica} was
not known; it was considered lost until 1896, when Richard Sch\"one recognized it as part of
a manuscript kept in a library in Istanbul.\footnote{This was first published in \cite{Heron1903}.
A modern study of this sole known surviving manuscript of the \emph{Metrica} can be found in 
\cite{Heron2014}.}

To help make some comparisons between various proofs, we begin with a version of the diagram in 
Heron's proof for a specific triangle.\footnote{See Figure 1.  To generate these figures, we used the triangle with vertices at $B = (0,0)$, $C = (5,0)$ and $A = (1,2)$ in the Cartesian plane.  This happens to have a right angle at $A$ so some of the line
segments in the figures are in rather special positions that facilitated the plotting.  However, this does not affect the arguments.  None of the authors we consider would have done things this way, of course.} 
\begin{figure}
  \centering
  \includegraphics[width=.6\linewidth]{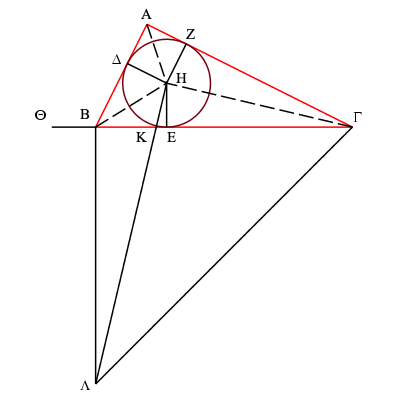}
  \caption{Heron's diagram in the \emph{Metrica}.}
\end{figure}
Heron's proof in outline consists of the following steps.  First, let the circle $Z\Delta E$ with center at $H$ 
be inscribed in the triangle $AB\Gamma$.  (Proposition 4 in Book IV of Euclid gives a construction for this where 
$H$ is found as the intersection of two of the angle bisectors of the triangle, but Heron takes this as known and
does not mention it explicitly.)   $H\Delta = HZ = HE$ since these are all radii of the same circle. Then since $HE$, $H\Delta$ and $HZ$ are perpendicular to the sides of the triangle, the area of triangle $AB\Gamma$ will equal
one half times $EH$ times the perimeter of the triangle.  Heron then does further construction steps, first extending 
$B\Gamma$ to $\Theta\Gamma$, letting $\Theta B = A\Delta$.  This makes $\Theta\Gamma$ equal to the 
semiperimeter of the triangle.   Second, he takes $H\Lambda$ perpendicular to 
$H\Gamma$ and extends to $\Lambda$ which is the intersection with the line through $B$ perpendicular to 
$B \Gamma$.  It follows that $HB\Lambda\Gamma$ is a cyclic quadrilateral and facts about the diagonals in
such quadrilaterals imply triangle $B\Lambda\Gamma$
is similar to triangle $\Delta HA$.  The proportionality of corresponding sides implies the square of the area of the triangle 
$AB\Gamma$ is equal to the square on the perpendicular $HE = H\Delta = HZ$ above times the square on 
the semiperimeter (using the fact that $B\Gamma$ and $\Theta B = A\Delta$ together equal the semiperimeter.)

Several of the previous practical geometry texts that Clavius mentions in his preface (see above) also 
include proofs of Heron's process/formula for finding the area of a triangle.  \emph{Significantly for 
our question, and reflecting the fact that the \emph{Metrica} was not known directly, none of them 
attributes this to Heron either}.   One of the earliest that does is
the groundbreaking \emph{De Practica Geometrie} of Leonardo of Pisa (``Fibonacci'') (ca.~1170--ca.~1250).   The
later \emph{Summa de arithmetica geometria proportioni et proportionalit\`a} by Luca Pacioli (1447--1517) does 
as well and Pacioli's treatment is virtually a copy of what Fibonacci says (although written in the Tuscan dialect
of Italian rather than Latin).\footnote{I have consulted \cite{Pacioli}, a scanned version of the 1523 edition of 
Pacioli's book at {\tt www.e-rara.ch}.}

  Marshall Clagett has written that Leonardo ``borrows heavily
and often in verbatim fashion'' in the revised version of this work from the \emph{Verba filiorum Moysi 
filii Sekir, i.e. Maumeti, Hameti, Hasen}.\footnote{See \cite[p.~224]{Clagett}.}
This work, also known as the \emph{Liber trium fratrum de geometria}, is a Latin translation of an Arabic work on mensuration by the 9th century Ban\=u M\=us\=a brothers made by Gerard of Cremona (1114--1187).  The original authors were key figures in the early translation movement by which Greek mathematics
became known in the Islamic world and the Greek original of the
\emph{Metrica} may have been available in Baghdad at this time.  But there are significant differences 
between the \emph{Metrica} version and the \emph{Verba filiorum} version.
Figure 2 shows what the diagram in the \emph{Verba filiorum} looks like for our triangle:
\begin{figure}
  \centering
  \includegraphics[width=.6\linewidth]{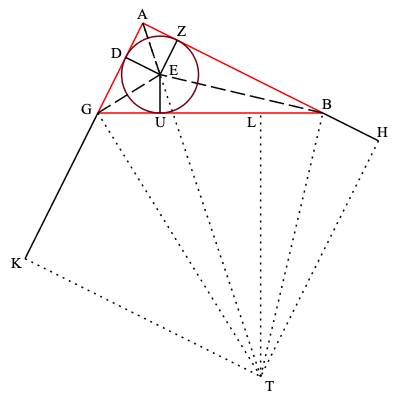}
  \caption{The \emph{Verba filiorum} diagram.}
\end{figure}
After identifying the points $D$, $Z$, $U$ as the points of tangency of the inscribed circle, the side $AB$
is extended to $AH$ by making $BH = GU$, so $AH$ is equal to the semiperimeter.  Similarly $AG$ is 
extended to $AK$ making $GK = BU = BZ$ and angle $AKT$ is a right angle.  The point $T$ is 
chosen so that it lies on the angle bisector at $A$.  It follows that triangle $EBU$ is similar to triangle
$BTH$ and the proportionality of corresponding sides implies the square of the area of the triangle 
$AB\Gamma$ is equal to the square on the perpendicular $EU = EZ$ times to the square on 
the semiperimeter $AH$.

Clagett  mentions that he believes Fibonacci's debt
to the Ban\=u M\=us\=a applies specifically to the treatment
of Heron's formula in Fibonacci's work.\footnote{See \cite[p.~224]{Clagett}.}  
However, a close analysis of the argument and the diagrams 
provided shows that while the proof of Heron's formula
in Proposition VII of the \emph{Verba filiorum}, and reproduced in Clagett's book, has many features in common with
Fibonacci's proof, it also has other features in common with the proof from Heron's \emph{Metrica} that do not occur in Fibonacci.  Probably the major example here is that the whole first phase of the argument in both the \emph{Metrica} proof and in the \emph{Verba filiorum} proof consists of considering the inscribed circle in the triangle (as in Proposition 4 
from Book IV of Euclid).  Fibonacci does not mention the inscribed circle; in fact he ends up repeating a large 
portion of the Euclidean proof to show that if perpendiculars (or as Fibonacci says, ``cathetes'') are dropped to the three sides from the intersection point of two angle bisectors in the triangle, then the three perpendicular segments are equal.\footnote{There is an unfortunate mistranslation at the start of the proof of Heron's formula in \cite{Fibonacci2008}.  At the start of the first full paragraph on p. 81,
the Hughes translation says, ``To prove this: in triangle $abg$ bisect the two equal angles $abg$ and $agb$ ... .''  This would make
the proof apply only to isosceles triangles.  But that is not correct.  The Latin text in the 19th century version
edited by B. Boncampagni, \cite[p. 40]{Fibonacci1862} says at this point:  ``Ad cuius rei demonstrationem adiaceat trigonum $abg$: et 
dividantur in duo equa anguli, qui sub .abg. et .agb. a rectis .bt. et .tg. ... ''  That is, ``To prove this: in the triangle 
$abg$, let the angles $abg$ and $agb$ [each] be divided into two equal angles by the lines $bt$ and $tg$ ... ''  Fibonacci
is definitely \emph{not} restricting his discussion to isosceles triangles.}
There are also some less drastic differences in the way that similar triangles within the figure are used to 
deduce that the square of the area is equal to the square on the perpendicular above times the square on 
the semiperimeter (and the square on the semiperimeter equals the semiperimeter times the product of the three excesses of the semiperimeter over the sides).
So it is surely not entirely accurate to characterize Fibonacci's proof (at least as a whole) as ``verbatim borrowing'' even if the overall strategies of the proofs are similar and the final sections of the proofs do more or less converge.\footnote{The difference in the diagrams is also mentioned in
\cite{Fibonacci2008}.  See the footnote on p. 83.}

On the other hand, Clavius's version of the proof of Heron's formula is different again, but significantly closer to the proof in Fibonacci than it is to the proof in the \emph{Verba filiorum}.  To discuss this in more detail it will be 
necessary to consider the  diagrams from these two proofs.  (See Figures 3 and 4.)   These two figures show Clavius's and Fibonacci's constructions
applied to the same particular triangle as in the previous figures.
\begin{figure}
\centering
\begin{minipage}{.5\textwidth}
  \centering
  \includegraphics[width=.95\linewidth]{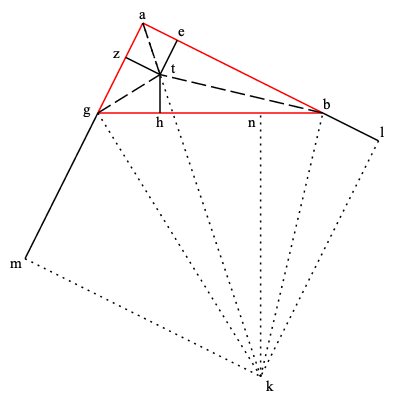}
  \captionof{figure}{Fibonacci's diagram.}
  \label{fig:test1}
\end{minipage}%
\begin{minipage}{.5\textwidth}
  \centering
  \includegraphics[width=.95\linewidth]{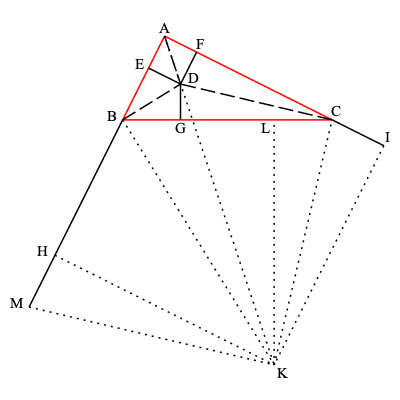}
  \captionof{figure}{Clavius's diagram.}
  \label{fig:test2}
\end{minipage}
\end{figure}  
Both Clavius and Fibonacci start by considering angle bisectors (these are two of the dashed black lines) 
for the two vertices on the horizontal side in the diagram and their intersection point 
($t$, and $D$, respectively).  They both
drop perpendiculars ($th$, $tz$, $te$, and $DG$, $DE$, $DF$ resp.) and use facts about congruent triangles
in the figure to show that the three perpendiculars are equal, that $at$ (resp. $AD$) also bisects that angle and, moreover, the two segments closest to each vertex are equal--that is $ae = az$, $be = bh$, $gz = gh$ (resp. $AE = EF$, $CF = CG$, and $BE = BG$).  
Neither mentions the inscribed circle, which would be tangent to the sides of the triangle in the points  
$e,z,h$ (resp. $E,F,G$). This 
implies that any one of the sides, together with one of the equal segments not meeting that side are together 
equal to the semiperimeter of the triangle -- for example, side $ab$ (resp. $AC$) together with $gh$ or $gz$ 
(resp. $BG$ or $BE$). In addition, the three excesses of the semiperimeter over the sides of the triangle 
that feature in Heron's formula coincide with the segments:  $ae$ or $az$, $be$ or $bh$, $gh$ or $gz$
(resp.  $AE$ or $AF$, $BG$ or $BE$, $CG$ or $CF$).  

Then, in further parallel constructions, the sides $ag$, $ab$ (resp. $AB$ and $AC$) are 
extended to $am$, $bl$ (resp. $BH$, $AI$) by making $gm = hb$ and $bl = gh$ (resp. $BH = GC$ and $CI = BG$).
As noted before, this makes both $am$ and $al$ (resp. $AH$ and $AI$) equal to the semiperimeter of the 
triangle, hence equal.\footnote{In the diagrams here, in fact, they are two sides of a square, but that is 
only true because our triangles have right angles at $a$, $A$.}  At this point, Fibonacci says to produce
the third angle bisector $at$ until it meets the segment $lk$ making a right angle with $ab$ at $k$.   Clavius, on the 
other hand (literally!), says to produce $AD$ to $K$ where it meets the line through $H$ perpendicular to $AH$.  
But either way, the next deduction is that by the SAS criterion, the triangles $amk$ and $alk$ 
(resp. $AHK$ and $AIK$) are congruent, so angle $amk$ (resp. $AIK$) is also a right angle, and moreover
$mk = lk$ (resp. $HK = IK$).  

In the final constructions, Fibonacci says to cut off the segment $bn$ from $gb$ so that 
$gn = gm = bh$, and hence $bn = bl = gh$.  Clavius does the parallel operations making 
$BL = BH = CG$ and hence $CL = CI = BG$.  But now Clavius does one further step that 
Fibonacci does not:  He extends $AH$ to $AM$, making $HM = CL = CI$. With $k$, resp. $K$ joined to all of
the newly constructed points, both proofs proceed to show that the lines $kn$ (resp. $KL$) meet
the horizontal side in a right angle.  The additional triangle $HMK$ introduced in Clavius's
argument is congruent to triangles $CIK$ and $CLK$ (resp. $bnk$, $blk$ in Fibonacci's figure), and hence it
is somewhat redundant.  But what we have here would seem to be a typical kind of procedure
for Clavius; at the cost of a few more steps, he furnishes a reader of his proof with another triangle 
$HKM$ that gives a perhaps easier way to understand why the angle at $n$ or $L$ is a right angle (this is
not really clear visually in Clavius's original diagram, where it seems no attempt has been made to 
show all the right angles accurately).  

Finally, as in all of the proofs of Heron's formula we have discussed, similar triangles can be identified in the figure
such that the proportionality of corresponding sides implies that the square of the perpendicular (e.g. $DE$ in Clavius's figure)
times the semiperimeter (e.g. $AH$ in Clavius's figure) is equal to the product of the 
three excesses of the semiperimeter over the sides, or $(DE)^2 (AH) = (EB)(BH)(AE)$.  
Hence the square of the perpendicular 
times the square of the semiperimeter, that is, $(DE)^2(AH)^2$, equals 
$(AH)(EB)(BH)(AE)$.  With some rearrangement of factors, this equals $s(s - a)(s - b)(s -c)$  in the modern
algebraic form of Heron's formula.  Fibonacci uses the triangles $ebt$ and $kbl$;
Clavius uses the triangles $AED$ and $AHK$.  The details in this step are somewhat different, but 
the idea is analogous.  

The final step in both of these proofs is to note that by the usual ``one half base times height'' way of computing 
areas of triangles, the sum of the areas of the three triangles $atb$, $btg$, $gta$ (resp. $ADC$, $CDB$, 
$BDA$), which equals the area of the whole triangle, is also equal to the product of any one of the three equal
perpendiculars and the semiperimeter.  Since Clavius has not proved the basic method for computing 
areas of triangles yet in the \emph{Geometria Practica}, 
he has to give a forward reference to the next section in Chapter 2 of his Book IV.  Fibonacci puts the
discussion of Heron's formula
after his treatment the other method, so he has set up what he needs already.

To conclude this section, we can say that Clavius did not attribute this result to Heron because
in his time it was simply not generally known that an ancient source for this result was the 
\emph{Metrica} of Heron.   On the other hand, 
let us offer the conjecture that here, as was often the case, Clavius reworked
and amplified what he found in other sources so that his version has additional or alternate features
intended to heighten clarity or to increase convenience for his readers.  Here it seems very probable that he was 
looking at Fibonacci's proof (or perhaps other proofs derived from that one, such as the 
proof in Pacioli's text) but his version is not a verbatim copy, any more than Fibonacci's was
a verbatim copy of the proof in the \emph{Verba filiorum}.  Why Clavius chose not to say this explicitly at this
point in his book is still somewhat mysterious, however.  In analogous situations, Clavius did sometimes
say explicitly how his account of a proof would differ from what was found in his source(s).\footnote{For example, 
he was very explicit about this in the introduction to Book VI on \emph{geodesy}.}

\section{Clavius's treatment of Archimedes' ``Measurement of the Circle'' in Book IV.}
\label{ArchimedesMeasurementCircle}

After Greek versions of this work of Archimedes (including summaries from Book V of the \emph{Mathematical 
Collection} of Pappus and the commentary on Ptolemy's \emph{Almagest} by Theon) 
were intensively studied in the Islamic world and the resulting Arabic translations were retranslated into Latin, the \emph{Measurement of the Circle} was surely 
the best-known and most-copied Archimedean text throughout the medieval period in western Europe. 
A major part of the reason for this was certainly the utility of the results of this work for practical questions.  
 On the other hand, 
the brevity of the work and its somewhat sketchy form have led Dijksterhuis to conjecture that ``it is quite
possible that the fragment
we possess formed part of a larger work,''\footnote{See \cite[p. 222]{Dijksterhuis}.}
and Knorr to judge that the versions we have represent 
 ``at best an extract from the original composition.''\footnote{See \cite[p. 375]{Knorr1989}.}
In \cite{Clagett}, Clagett reproduces two translations of this work from Arabic into Latin, the first made (``perhaps'') by 
Plato of Tivoli (fl.~12th century), and the second made by the same Gerard of Cremona mentioned in the 
previous section.  Clagett also reproduces six additional ``emended'' versions as well as the 
treatment of the results of this work in the \emph{Verba Filiorum}, following the Ban\=u M\=us\=a.  
Part III of Knorr, \cite{Knorr1989}, contains a more complete study of the transmission including additional
versions and reflecting more recent scholarship.   Here 
we will simply say that the common elements of most of these are three propositions stated in this order and  in something like these forms:

\begin{prop}[Archimedes, \emph{Measurement of the Circle}, 1]
Every circle is equal [in area] to a right triangle, one of whose sides containing the right angle is equal to the radius
of the circle while the other side containing the right angle is equal to the circumference.
\end{prop}

\begin{prop}[Archimedes, \emph{Measurement of the Circle}, 2]
The ratio of the area of any circle to the square of its diameter is the ratio $11$ to $14$.
\end{prop}

\begin{prop}[Archimedes, \emph{Measurement of the Circle}, 3]
The circumference of a circle exceeds three times its diameter by a quantity less than $\frac{1}{7}$ of the diameter and 
greater than $\frac{10}{71}$ of the diameter.\footnote{That is, in modern terms, $3\,\frac{10}{71} < \pi < 3\,\frac{1}{7}$.  
Archimedes may well have used methods similar to the ones to be discussed to produce tighter estimates
for the ratio of the circumference to the diameter.  But if so, no text doing this has survived.}
\end{prop}

It is interesting to compare Clavius's treatment of these results to what authors of other 
practical geometry texts say.  One of the earliest, the \emph{Practica Geometriae} attributed to Hugh of 
St. Victor (ca.~1096--1141), \cite{Hugh}, does not mention this topic at all.  
In his \emph{De Practica Geometrie}, Fibonacci addresses 
the content of all three propositions in turn.\footnote{See \cite[pp. 152--158]{Fibonacci2008}, paragraphs [191]-[200] of Chapter 3.}  
However, Fibonacci does not really attempt to present the full  
proof for Proposition 1 that is found in other versions.
Instead, what he says is as follows.  First, 
considering a regular polygon circumscribed about the circle, Fibonacci argues that the product
of the radius and the perimeter of the polygon is greater than the area of the circle
by considering the triangles formed by joining the center and the vertices of the polygon.  
\begin{figure}[h]
 \centering
  \includegraphics[width=.5\linewidth]{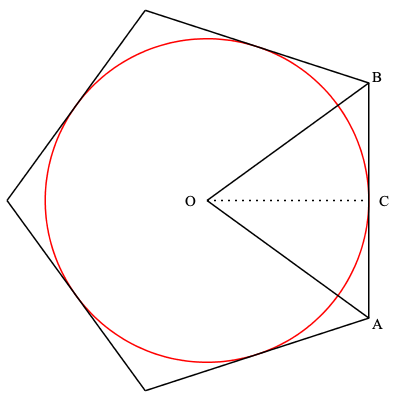}
  \caption{One of Fibonacci's circumscribed polygons.}
\end{figure}
For example, the area of the pentagon in Figure 5 will be $5 \cdot \frac{1}{2} \cdot (AB)(OC)$, which is greater
than the area of the circle.  Hence
the product of the radius and a number greater than half the circumference of the circle gives
an area greater than the circle.  Fibonacci apparently
takes it as obvious that the perimeter of the circumscribed polygon is greater than the circumference
of the circle.\footnote{Clavius returns to this point in Book VIII of the \emph{Practical Geometry}, discussing
arguments by Archimedes from \emph{On the Sphere and Cylinder}, and an alternate treatment by Girolamo
Cardano.} 

Next, by a clever observation, Fibonacci
considers an inscribed $n$-gon and adds vertices bisecting the arcs between successive vertices
of the $n$-gon to form an inscribed regular $2n$-gon. 
\begin{figure}[h]
 \centering
  \includegraphics[width=.5\linewidth]{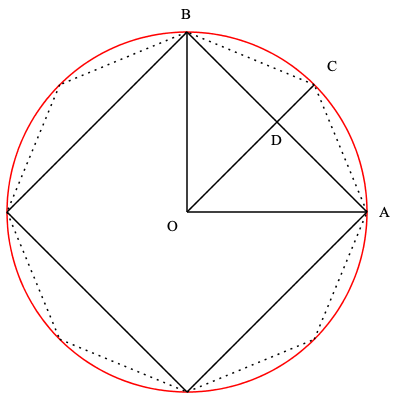}
  \caption{Fibonacci's inscribed polygons with $n = 4$ and $2n = 8$.}
\end{figure}
 Fibonacci notes that the product of the radius and  
half the perimeter of the $n$-gon is equal to the area of the $2n$-gon, hence less than the area of the circle.  
This follows since, for instance, the product $\frac{1}{2} (OC)(AB)$ in Figure 6 is equal to the sum of the areas
of the triangles $OAB$ and $ABC$.  Hence four times this will equal the area of the octagon drawn with dotted lines.
Similar arguments apply for all $n\ge 4$.
Therefore multiplying the radius by a number less than half the circumference of the circle gives an
area less than the area of the circle.  ``Whence it is concluded that the product of the radius of the circle
and half its circumference equals its area.''\footnote{Quare concluditur, quod ex multiplicatione semidiametrij
circulj in dimidium lineae circumferentis prouenit embadum ispius. See \cite[p. 87]{Fibonacci1862}.}  

Although it is certainly intuitively clear that the areas of the inscribed and circumscribed polygons converge to the 
area of the circle as $n\to \infty$, it seems fair to characterize what Fibonacci says as more of a plausibility argument 
than a complete proof because he has not shown that the difference between a circumscribed polygon and an inscribed
polygon can be made arbitrarily small.  

Second, Fibonacci argues by way of Euclid XII.2 that the ratio of the square of the diameter to the area is the 
same for all circles.  He essentially then does a ``proof by numerical example''\footnote{fuit enim quadratum 
dyametrij suprascripti $196.$ et embadum ipsius $154.$ quorum proportio est sicut $14.$ ad $11.$ ...  See 
\cite[p. 88]{Fibonacci1862}.} for Proposition 2, using the result of Proposition 3 and effectively taking $\pi = 22/7$.   There is 
no indication that this ratio is only an approximation and that no actual circle has the square of the diameter
exactly equal to $196$ and area exactly $154$.  

Finally, Fibonacci turns to the Archimedean estimates $3\,\frac{10}{71} < \pi < 3\,\frac{1}{7}$ (and here mentions 
Archimedes explicitly for the first time).  He says that he is not going to follow Archimedes' proof exactly because
smaller numbers will suffice to make the point.\footnote{Ostendendum est etiam quomodo inuentum fuit lineam
circumferentem omnis circulij esse triplum et septima sui dyametrij ab Archimede philosopho, et fuit illa inuentio
pulcra et subtilis ualde: quam etiam reiterabo non cum suis numeris, quibus ipse usus fuit demonstrare; cum
possibile sit cum paruis numeris ea que ipse cum magnis ostendit plenissime demonstrare.  See
\cite[p. 88]{Fibonacci1862}.  Note that again there is no indication that $22/7$ is only an approximation
to the ratio in question.}  Interestingly, Fibonacci includes more about the individual details of the calculations
paralleling the proof of Proposition 3 than most of the versions of the \emph{Measurement of the Circle} do.  

As was true for the discussion of Heron's formula considered in the previous section, Luca Pacioli's 
discussion of the results of the \emph{Measurement of the Circle} follows Fibonacci very closely, even with
virtually identical diagrams.\footnote{For instance, about the proof of Proposition 3, Pacioli says 
``Ancora eglie da mostrare comme
e so trouata da Archimenide la linea circonferentiale essere $.3.$ volte $.\frac{1}{7}.$ del diametro: la
quale inventione so bella e sotile in questo modo, bene che con breuita se dica,'' \cite{Pacioli}, Pars Secunda, Distinctio quarta, Capitulum secundum, folio 31. (In this printed book, one folio, i.e. 
two facing pages, \emph{verso} and \emph{recto}, are given just one page number).}
Other practical geometry texts after Clavius's time tended to include the Archimedean estimates and discuss how to compute areas of circles and parts of circles, but to omit proofs for these facts entirely.

By contrast with Fibonacci or Pacioli, Clavius 
explains what he aims to do in this introductory paragraph:
\begin{quote}
It will not be a digression, therefore, if I include [Archimedes'] truly 
most acute and precise book, partly
because it is very brief (indeed, it consists of only three propositions), 
partly so that the student, in order to understand something so useful and
so widely applied in the works of all authors, should not be forced to go to 
Archimedes himself, and finally mostly because the writings of Archimedes,
as a result of their brevity, are somewhat obscure, and we hope to 
bring some light to them.\footnote{Non abs re ergo sit, si eius libellum
de circuli dimensione acutissimum sane, \& subtilissimum hic interferam, tum quia breuissimus
est, quippe qui tribus duntaxat propositionibus constet: tum ne studiosus, vt rem tam vtilem,
atque apud omnes artifices peruulgatam intelligat, Archimedem ipsum adire cogatur: tum 
vero maxim\`e, quod cum Archimedis scripta ob affectam breuitatem, sint paulo obscuriora,
illis nos lucem aliquam allaturos speramus.}
  (\cite[p.~182]{Clavius1606})
\end{quote}

Clavius starts with an account of the rather subtle exhaustion proof of Proposition 1 
found in most of the versions of the \emph{Measurement of the Circle} mentioned above. It must 
be proved that the area of a circle is the same as the area of a right triangle with the two sides about the 
right angle equal to the radius of the circle and the circumference of the circle.\footnote{This is very closely related
to the results on isoperimetric problems to be discussed later in Book VII.}   This proof has much in
common with the proof of Proposition 2 in Book XII of Euclid.  The plan is to show
that assuming the area of the circle is either greater than or less than the area of a right triangle as in the statement 
leads to a contradiction.  A key role will be played by a statement introduced by Euclid in the proof of 
Proposition 1 in Book X of the \emph{Elements}.
Clavius uses this in the form: \emph{If an area at least half the area of the
circle is taken away from the circle, and from the residual
area again an area at least half of that remaining area is
taken away, and so on, a there will eventually remain an
area less than any positive magnitude $z$.}

First suppose the circle is larger than the stated triangle
by a certain positive magnitude (Clavius calls this $z$).
Let a sequence of non-overlapping areas be removed from the interior 
of the circle.   Specifically, the inscribed square is removed at the first step, then
 four triangles on the sides of the square inscribed in the circle, so at this point a
regular octagon has been removed, then eight triangles on the sides
of the octagon, so the total removed figure is a regular 16-gon, etc.  
Clavius actually stops with the octagon so he has a diagram that is 
virtually identical with the ones presented by Fibonacci and Pacioli.\footnote{See Figure 7.  Many of 
  Clavius's point labels and additional lines constructed in Clavius's very ``busy'' diagram have been omitted at this stage for clarity.}
\begin{figure}[h]
 \centering
  \includegraphics[width=.45\linewidth]{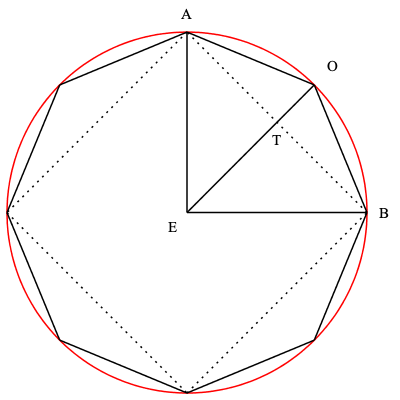}
  \caption{The essential portions of Clavius's diagram for the proof of the first part of Proposition 1.}
\end{figure}
Note that the right triangle $EAB$ is half of the square with $E,A,B$ as three vertices. Four of these triangles
make up the inscribed square in the circle and four of those small squares make up a square circumscribed 
about the circle.  So the inscribed square is half of the circumscribed square and hence more than half
of the circle.  Similarly, the four triangles on sides of the inscribed square (four triangles congruent to $ABO$)
are half the rectangles with bases equal to the sides of the inscribed square and heights equal to the perpendicular segments such as $OT$.  Hence removing those four triangles removes more than half of the area between
the inscribed square and the circle.  It is not hard to show that this pattern continues indefinitely.  So 
eventually the remaining region between a $2^m$-gon and the circle will have area less than the positive
magnitude $z$.  However, this leads to a contradiction if the construction is applied sufficiently many times.
The area of the circle is supposed to equal $({\rm area\ of\ triangle}) + z$, but
the area of the circle also equals:
$$({\rm area\ of\ inscribed\ polygon}) + ({\rm remaining\ area}) < ({\rm area\ of\ triangle}) + z,$$
since the perimeter of the polygon is smaller than the circumference of the circle,
and the apothem--the perpendicular from the center to the side of the inscribed polygon--is less than the radius, 
and hence the area of the inscribed polygon is less than the area of the triangle. 
Hence the area of the circle cannot be greater than the area of the triangle.

Now suppose that the circle is less than the stated triangle by 
a certain magnitude (again denoted $z$).
Clavius and Archimedes now start from the square circumscribed  about
the circle (whose area is definitely greater than the area of the triangle since the perimeter
of the square is greater than the circumference of the circle, and the 
apothem is the same as the radius).
Clavius begins removing areas from the square: first the circle, then four exterior
triangles (such as $KXV$ in Figure 8) with base tangent to the circle, then again 
eight exterior triangles with base tangent to the circle at the midpoints of the arcs $BO$, $OV$, and so forth.  
The remaining regions (after the very first step when the circle is removed) now
are collections of what Clavius calls ``mixed triangles,'' with one side an arc of the circle.
\begin{figure}[h]
 \centering
  \includegraphics[width=.45\linewidth]{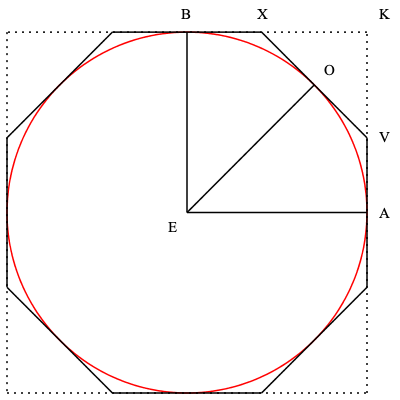}
  \caption{The essential portions of Clavius's diagram for the proof of the second part of Proposition 1.}
\end{figure}
Again since the inscribed square in the circle is half of the circumscribed circle, after the circle 
has been removed, less than half of the circumscribed square remains.  Similarly, each of the four 
exterior triangles such as $KXV$ is greater than half of the area between the circle and the lines
$BK$, $KA$, and so forth.   Since more than half the remaining area is removed at each step, 
eventually the remaining area becomes less than $z$.
But this also leads to a contradiction.  On the one hand
$$({\rm area\ of\ circumscribed\ polygon}) > ({\rm area\ of\ triangle})$$
because the perimeter of the polygon is greater than the circumference and the apothems are now
all equal to the radius of the circle.  But on the other hand, by the process described above,
$$({\rm area\ of\ circumscribed\ polygon}) < ({\rm area\ of\ circle})+ z = ({\rm area\ of\ triangle}).$$
Hence, since the triangle in the statement of the proposition is neither greater nor less than
the circle, it can only {\it equal the circle}. 
Clavius's version of this proof incorporates many
explanatory comments and justifications for the individual steps in the reasoning not found in other
versions.  

An amusing sidelight in the form of a long Scholium
follows, in which Clavius refutes Joseph Justus Scaliger's\footnote{Scaliger (1540--1609) was an eminent French Protestant classical philologist and historian who 
also fancied himself a mathematician.  He managed to convince himself both
that the area of a circle is $\frac{6}{5}$ times the area of the inscribed regular hexagon, a statement 
that is equivalent to $\pi = \frac{9\sqrt{3}}{5}$, \emph{and also} that the area of a square 
with side equal to the circumference of a unit circle is ten times the area of the square on the diameter, a statement equivalent to $\pi = \sqrt{10}$.     
These are not consistent with each other and neither is consistent with the Archimedean estimates
$3\,\frac{10}{71} < \pi < 3\,\frac{1}{7}$. Moreover, either one
would imply that the circle can be squared with straightedge and compass, which we know today is 
impossible.  Clavius suggests these misunderstandings as one possible reason
for Scaliger's claim that the Archimedean argument
cannot be correct. Scaliger's faulty conclusions were
included in his mathematical  
{\it magnum opus}, grandly titled {\it Cyclometrica Elementa}, published in 1594 in a lavish edition
with statements of theorems in both Latin and ancient Greek.  In 1609 Clavius
published an 84-page pamphlet \emph{Refutatio Cyclometriae Iosephi Scaligeri} (Refutation of the \emph{Cyclometrica} of Joseph Scaliger), giving a blow-by-blow analysis of all of the (numerous) errors in Scaliger's work.
This is contained as an appendix in Volume V of the \emph{Opera Mathematica}. } claim that Archimedes must have been mistaken
in this proof.  Indeed, Clavius takes Scaliger to task rather savagely over his misunderstandings.  A small sample:

\begin{quote}
 And I am really astonished that you,
\emph{Mathematicus} that you are, deny that some quantity is equal
to another when it is neither greater nor less.  For if it is not equal to the 
other, then it will be unequal to the other, therefore either greater or 
less, against that hypothesis.  Or don't you see that not only Archimedes
but also Euclid used this way of arguing most frequently in Book XII of the \emph{Elements}?\footnote{Et 
sane miror, te, Mathematicus, cum
sis, negare quantitatem aliquam illi esse aequalem, qua neque maior est, neque minor.
Si enim aequalis non est, erit inaequalis.  Igitur vel maior vel minor, contra hypothesim
cum dicatur neque maior esse, neque minor.  An non vides, non solum Archimedem, sed
etiam Euclidem lib.~12.~hunc argumentandi modum frequentissim\`e vsurpare?}
(\cite[p.~185]{Clavius1606})
\end{quote}
Scaliger had
had another ``run-in'' with Clavius over the Gregorian calendar reform in which Clavius had taken a leading role, as well as ongoing controversies with other Jesuits on various subjects, so there was ample bad blood between them.  Some of that is manifest in the scathing polemical tone of Clavius's comments.  

  Following this, Clavius notes that the usual Proposition 3 is used in the proof of Proposition 2, and hence 
he has decided to reverse the order of Propositions 2 and 3 as found in other versions to maintain the chain 
of logical implications.\footnote{Clavius 
writes:  Haec est Archimedis propositio 3.~quam nos secundam facimus, vt doctrinae ordo 
servetur, quando quidem sequens propositio 3.~quam ipse 2.~facit, hanc nostram
propositionem 2.~in demonstrationem adhibet.}  In his account of this famous proof, Clavius essentially 
follows the plan used in most other versions.  There are again two steps, the first (in effect) considering 
polygons circumscribed about the circle and triangles
with one side along a tangent to the circle, the second (in effect) considering polygons inscribed in the 
circle and triangles inscribed 
a semicircle.  In each phase, an angle is repeatedly bisected, until one side in the triangle
comes from a regular $96$-gon.  The Pythagorean theorem is applied repeatedly to estimate ratios and
the proof consists essentially of a complicated series of numerical calculations, quite different from 
many Greek geometric arguments.

Two aspects of Clavius's version of the proof of his Proposition 2 ($=$ usual Proposition 3),
as compared with other versions, are notable.  First,
Clavius provides \emph{more details}, more \emph{justification for individual steps}, and 
\emph{a fuller treatment of the calculations} leading to the 
estimates than Archimedes (or whoever wrote the versions of the Archimedean text that we have) did.
Since the bisection steps in each half of the proof follow exactly the same plan, some versions
work out the first step in detail, and then just present the numerical results for the subsequent 
steps.\footnote{The author used this expedient, in fact, in his translation of Clavius's proof because 
the repetitive structure of the original makes it very tedious reading!}
This is very much along the lines of how he treats information from his sources in other sections of the 
\emph{Geometria Practica}.  The second aspect is probably somewhat less mathematically 
significant, but still interesting.  Namely, 
by showing the bisections from the first phase of the proof to the left of a vertical diameter
in the circle, and the bisections from the second phase to the right, Clavius manages to condense 
all the steps of the constructions for both phases of the proof into a single diagram in a clever, but
still entirely understandable way.\footnote{See the figure inserted for the first time on page [186] of 
the \emph{Geometria Practica} and repeated several times thereafter for the convenience of 
the reader.  Setting up the figure this way
could have been a purely practical decision intended to reduce the number of different figures that had to be produced in
printing the book, of course.}  Most versions provide separate diagrams for 
each phase.\footnote{See for instance the facsimiles in \cite[p.~460, 463]{Knorr1989}.}

Another quite significant deviation from sources such as Fibonacci is that in 
his Proposition 3 ($=$ usual Proposition 2), Clavius explicitly adds the qualifier
\emph{approximately} (proxim\`e in the Latin) to the usual statement that the ratio of a circle to the square on the diameter 
is the ratio $11$ to $14$.\footnote{Circulus quilibet ad quadratum diametri proportionem habet, quam ad $11.$ ad $14.$
proxim\`e.}  As we noted above, Fibonacci did not do this in his version of this part of the Archimedean text.
Clavius's proof is essentially the same, though, using the approximation $\pi \doteq 22/7$.

Finally, it is interesting to note that after his account of Archimedes' results, Clavius also includes
closer approximations to $\pi$ later in Book IV, quoting results of his contemporary
Ludolph van Ceulen (1540--1610), and his student, Jesuit colleague, and successor as 
professor of Mathematics at the \emph{Collegio Romano}, Christoph Grienberger (1561--1636).  
Clavius states the equivalent of the bounds 
$$3\,\frac{14159265358979323846}{100000000000000000000} < \pi < 3\,\frac{14159265358979323847}{100000000000000000000}$$
and tries to give a practical ``spin'' for how these might be useful.  If one of these estimates is used
(for instance, the upper bound to parallel the $22/7$ value), then 

\begin{quote}
... the area of the circle will differ less from the true value than the area found from the 
Archimedean ratio.  But since it is more difficult to compute with large
numbers than with small ones, in practice the Archimedean ratio is
applied.  However, when more accurate values are desired, the Ludolphine
ratio above should be used, especially for large circles.\footnote{... quae 
quidem area minus \`a vera distabit, quam illa, quae ex proportione Archimedis inuenitur.  Sed quia
difficilius est per magnos numeros calculum instituere, quam per minores, vsus artificum obtinuit,
vt proportio Archimedis ad calculum adhibeatur.  Quando tamen desideratur accuratior calculus,
vtendum erit posteriori hac proportione Ludlophi, praesertim in maioribus circulis.}
(\cite[p.~199]{Clavius1606})
\end{quote}
Even today, though, it is difficult to imagine where 20 decimal place accuracy might really be needed(!)

\section{Clavius's discussion of methods for finding two mean proportionals between two given lines in Book VI.}

In this section and the next, we will discuss two connected topics from Book VI of the \emph{Geometrica Practica}.  
The first is the extended solution of the following Proposition 15 from pages 266--272:

\begin{prop}
To find two mean proportionals between two given lines approximately.\footnote{Inter datas duas rectas, 
duas medias proportionales prope verum inuenire.}
\end{prop}

This, the closely connected problem of duplicating the cube, plus
the problems of squaring the circle and trisecting an arbitrary angle (which are discussed by Clavius 
in Books VII and VIII), 
were tremendous stimuli to Greek geometry over hundreds of years.\footnote{For an extensive modern study
of the surviving sources and the historical development, see \cite{Knorr2012}.}
As Clavius says,
\begin{quote}
We will first report what the
ancient geometers have left to us in their writings concerning this problem.  
For this drove and tormented the talents of many, although
up to this day, no one will truly and geometrically have constructed two
mean proportionals between two given lines.\footnote{Quocirca prius in hac
propos.~in medium afferemus, quae antiqui Geometrae nobis hac de rescripta
relinquerunt.  Multorum enim ingenia res haec exercuit, atque torsit, quamuis
nemo ad hanc vsque diem, ver\`e, ac Geometric\`e duas medias proportionales
inter duas rectas datas inuenerit.}  (\cite[p.~266]{Clavius1606}) 
\end{quote}

Hippocrates of Chios (ca.~470--ca.~410 BCE)
was traditionally credited with the reduction of the problem of duplicating the cube to the problem of 
constructing two mean proportionals between given lines.  If $AB$ and $CD$ are the lines, 
two other lines $XY$ and $ZW$ are said to be two mean proportionals (in continued proportion) if 
$$AB : XY :: XY : ZW \quad \text{and} \quad XY : ZW :: ZW : CD.$$
Representing the lengths by numbers and using algebra, this becomes the string of equations
$$\frac{AB}{XY} = \frac{XY}{ZW} = \frac{ZW}{CD},$$
from which it follows that 
$$\left(\frac{ZW}{CD}\right)^3 = \frac{AB}{CD}.$$
So for instance if  $CD = 1$ and $AB = 2$ in some units, a construction of the two mean proportionals gives the
line $ZW$ which has length $\sqrt[3]{2}$, and that is the edge length of the cube with twice the volume of the cube with
edge length $CD = 1$.   

Clavius has included what might seem to be a surprising amount of the Greek work on 
this construction problem in the 
\emph{Geometria Practica}.  Although 
the works of the authors involved did not survive from antiquity in their original forms, 
they were summarized and hence preserved in the commentary on Archimedes' \emph{On the Sphere and 
Cylinder} by Eutocius of Ascalon (ca.~480--ca.~540 CE).\footnote{This is translated together
with the Archimedean text in \cite{Netz}.}  This 
is Clavius's stated primary source for this material although he may also have consulted Book III of
the \emph{Mathematical Collection} of Pappus of Alexandria (ca.~290--ca.~350 CE) where some of the same
methods are surveyed.  

Note that Clavius explicitly says ``approximately'' in the statement of the problem 
(in the Latin:  \emph{prope verum}, literally ``near the truth'').
This feature might seem curious for modern readers and it surely 
deserves some elaboration.  Clavius makes this qualification because the 
solutions he will present all involve either limiting operations relying on the senses of the 
geometer (so-called \emph{neusis} (\textgreek{ne\~usis}) constructions) or the use of auxiliary curves such as 
\emph{cissoid} of Diocles or the \emph{conchoid}
of Nicomedes that cannot be drawn as a whole using only the straightedge and compass.\footnote{The same is 
true for the \emph{quadratrix} curve of Hippias that Clavius studies extensively in Book VII.}
Because these solutions use more than 
the traditional Euclidean tools, they don't qualify as what Clavius means by ``geometric'' or exact solutions.\footnote{In 
his well-known methdological discussions of different solutions of construction problems from Books III, IV, and VII of
the \emph{Mathematical Collection}, Pappus would say they are not \emph{planar} solutions. }  

It is understood today that no such purely ``geometric'' solutions are possible for the three problems 
mentioned above and it is primarily this \emph{methodological question}--are the three problems
solvable under the most severe restriction to the use of only the Euclidean tools?--that has 
survived in modern discussions.  For instance, many undergraduate algebra 
courses discuss these problems via coordinate geometry and the characterization of points constructible
with straightedge and  compass as those whose coordinates lie in a field at the end of a tower of 
quadratic extensions starting with the rational numbers. The fact that the index 
$[{\mathbb Q}(\sqrt[3]{2}) : \mathbb{Q}]$ is equal to $3$
shows that it is not possible to duplicate the cube with straightedge and compass, and that is often the end of the story.

But this is a very modern and ``pure mathematical''
way of looking at things.  For Clavius, as for at least some of the Greeks before him, 
although the methodological question might be interesting, it was also important to find \emph{some reasonably 
accurate method} for constructing the two mean proportionals even if it meant using approximate methods rather
than an exact, ``geometrical'' solution. Perhaps surprisingly, this is actually a very \emph{practical} problem that had important
applications in architecture, military science and many of the other areas Clavius mentions in his 
Preface.  It gives a method for determining the linear dimensions 
of a solid figure similar to a given figure whose volume has a given ratio to the volume of the given figure.  
Just as a procedure for finding one mean proportional lets one rescale a plane figure in a given ratio, a solution
for this problem lets one rescale solid figures in any given ratio, and Clavius points this out explicitly 
a number of pages later, after Proposition 17 in the same Book VI:
\begin{quote}
This establishes the method by which a cube is not only to be duplicated (which the ancients
were seeking), but also increased or decreased in any given ratio.  It also gives the 
method by which bores of cannons are to be made larger or smaller according to a given
ratio.\footnote{Constat ex his, qua ratione Cubus non solum duplicandus sit (quod veteres inquirebant) sed
etiam augendus minuendusue in quacunque proportione: Item quo pacto pylae bombardarum
maiores, aut minores fieri debeant secundum proportionem datam.  In this connection we also point out the first part of the
heading of the first method Clavius presents--Method of Heron in the introduction to the Mechanics
and Making of Missile-throwing Machines(!)}  (\cite[p.~274]{Clavius1606})
\end{quote}
We note that Fibonacci also discusses methods for finding two mean proportionals
in his \emph{De Practica Geometrie}.\footnote{See \cite[Chapter 5]{Fibonacci2008}, paragraphs [12]-[15].} 
We will return to this point shortly and compare his approach with Clavius's approach.

In introducing his discussion, Clavius says he is making a very deliberate choice from among the many solutions
presented in Eutocius's commentary:
\begin{quote}
Although they are most elegant and
acute, the solutions of Eratosthenes, Plato, Pappus of Alexandria, Sporus, 
Menaechmus by means of the hyperbola and parabola, then with the 
help of two parabolas, and Archytas of Tarentum will be omitted and 
we will explain only the four solutions from Heron and Apollonius of Perga, 
Philo of Byzantium and Philoponus, Diocles, and Nicomedes.  \emph{We have 
judged these to be more useful, easier, and less prone to error.}  Anyone who
should want the other methods will be able to read them in the commentary of 
Eutocius of Ascalon in the second book of \emph{On the Sphere and Cylinder} of 
Archimedes, and in the book of Johannes Werner of Nuremburg\footnote{German mathematician, 
1468--1522.} on the conic sections.\footnote{Praetermissis autem modis Eratosthenis; 
Platonis; Pappi Alexandrini; Spori; menechmi tum beneficio Hyperbolae, ac
parabolae, tum ope duarum parabolarum; \& Architae Tarentini, quamuis
acutissimis, subtilissimisque; solum quatuor ab Herone, Apollonio Pergaeo,
Philone Bysantio, Philoppono, Diocle, \& Nicomede traditos explicabimus,
quos commodiores, facilioresque, \& errori minus obnoxios iudicauimus.
Qui aliorum rationes desiderat, legere eas poterit in Commentarijs Eutocij
Ascalonitae in librum 2.~Archimedis de Sphaera, \& Cylindro: Item in libello
Ioannis Verneri Norimbergensis de sectionibus Conicis.} (\cite[p.~266]{Clavius1606}; emphasis added)
\end{quote}
In other words, the methods discussed here are sufficient for the applications Clavius has in mind and they
are the ones he thinks are easiest and best suited for practical implementation.

By way of contrast, Fibonacci makes a different 
selection and presents only the methods ascribed to Archytas, Philo, and Plato by Eutocius.  Hence there is very little
overlap between his account and Clavius's account.  Moreover, he presents the method of Archytas (which relies on 
some quite involved solid geometry) first, after saying that finding the two mean proportionals ``... is not a thing that can 
be done easily, but this is 
how it must be done.''\footnote{... hoc facili operari non possit, tamen, qualiter hoc fieri debeat.  \cite[p. 153]{Fibonacci1862}.}

Turning now to the details of Clavius's account, the first method presented actually combines
two very closely related approaches, ascribed to Heron and Apollonius and discussed separately
by Eutocius. 
\begin{figure}[h]
 \centering
  \includegraphics[width=.45\linewidth]{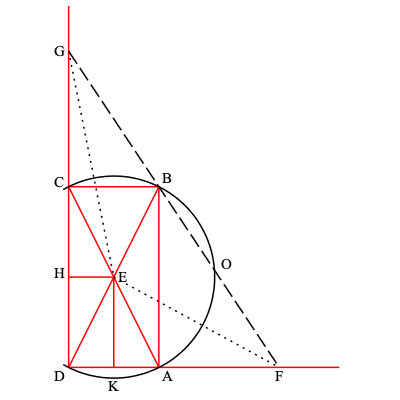}
  \caption{The essential portions of Clavius's diagram for the first two methods.  For Heron's method
  and Philo's method,  he dashed line $GBOF$ would need to be rotated about $B$ to reach the final desired position with 
  $EG = EF$ or $OF = BG$.}
\end{figure}
Clavius's version is a very close copy of Eutocius's text for Heron's method, with the variation 
represented by Apollonius's method inserted at one point.   In Figure~9, suppose we wish to 
find two mean proportionals between the lines $AB$ and $BC$, which have been arranged 
as two sides of the rectangle $ABCD$.  For Heron's method, Clavius says 
\begin{quote}
With sides $DA$, $DC$ extended, it is understood that a straightedge [represented by the dashed line
in the figure] placed 
at $B$ should be moved until it meets $DA$, $DC$, produced, in points $F$ and $G$
such that the lines $EF$ and $EG$ are equal.\footnote{Protractis autem lateribus, 
$DA$, $DC$, intelligatur circa punctum $B$, moueri
regula hinc inde, donec ita secet $DA$, $DC$ productas in $F$, \& $G$, vt rectae 
emissae $EF$, $EG$, aequales sint.} (\cite[p.~267]{Clavius1606})
\end{quote}
When this is true, consideration of the various similar triangles in the figure shows that 
$AF$ and $CG$ are the two desired mean proportionals between $AB$ and $BC$.  
Apollonius's variation of this method consists of finding a circle with
center at $E$ which has a chord $GF$ passing through $B$, and hence $EF = EG$ again.  
Clavius includes a brief description of a trial-and-error
method for finding the required circle not found in Eutocius.   

The second method, ascribed to Philo and Philoponus, has been reworked and greatly simplified by Clavius based on the realization that it is again very closely related to the first one (in fact Clavius sets up the discussion so that the same diagram applies).  Namely, with the circle $DABC$ described with 
center $E$ and radius $EA = EB = EC = ED$, the ruler at $B$ (that is, the dashed line in the figure)
is moved until $BG = OF$, where $O$ is the 
second intersection with the circle above.  Then it is easy to see we are back in exactly the same configuration
as in the other methods, so the same reasoning applies to give the two mean proportionals.   In our opinion,
this family of methods would certainly be among the easiest to apply.  They would probably be the most accurate
as well.  Given a sufficiently accurate diagram
of the rectangle $ABCD$ and its diagonals, rotating
a ruler passing through $B$, or using a compass to draw various test circles with center at $E$ would certainly give satisfactory results
if they were performed with sufficient care: the chances of making large errors would be extremely small since 
it is so clear what is required.  

Note that the geometer is required to rotate the line through
$B$ or adjust the radius of a circle centered at $E$ until a certain condition is satisfied.  As presented by Clavius,
this involves approximation processes making use of the senses of the geometer, as we said earlier.  The 
next two methods will be somewhat different in that they are set up to make use of \emph{auxiliary curves} 
whose description (that is, the description of the whole curve and not just a finite set of points on the curve)
requires tools besides the straightedge and compass.  
  
The next method Clavius discusses is ascribed Eutocius to Diocles (ca.~240--ca.~180 BCE), and specifically to a 
book called \emph{On Burning Mirrors}.  The Greek original has not survived so this was known only 
from fragments preserved in other texts like Eutocius's commentary.  But an Arabic translation 
of the whole has survived and this has now been translated into English by G.~J.~Toomer, \cite{Diocles}.  Clavius covers essentially the same ground as in the corresponding
section from Eutocius's commentary.  However, as usual, he has reworked and augmented his source material significantly.
Clavius begins by separating off what he calls the ``Lemma of Diocles,'' which identifies a geometric configuration
containing two mean proportionals between given lines.  
\begin{figure}[h]
 \centering
  \includegraphics[width=.45\linewidth]{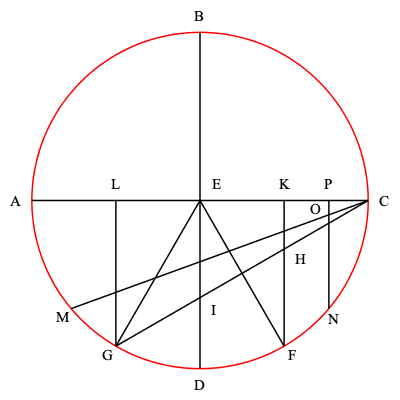}
  \caption{Clavius's figure for the ``Lemma of Diocles."}
\end{figure}
(See Figure 10. An equivalent figure with Greek letter labels appears in Eutocius.)
Let $AC$ and $BD$ be diameters of the circle meeting at right angles at $E$.  Let arcs $DG$ and $DF$ be 
equal and join $CG$.  Let $GL$ and $KF$ be drawn parallel to $BD$.  Let $CG$ meet $FK$ at $H$.  Then 
by considering relationships of the lines in the figure, Clavius essentially follows the proof given by Eutocius to show
that $FK$ and $KC$ are two mean proportionals
between $AK$ and $KH$.  Similarly, if the arcs $DM$ and $DN$ are equal, then drawing $CM$ cutting the vertical
line $PN$ in $O$, it follows that $NP$ and $PC$ are two mean proportionals between $AP$ and $PO$.

Now given two lines $AB > BC$, we can apply the ``Lemma of Diocles" in the following way:  First construct
a circle with radius $AB$ and lay off $BC$ along a perpendicular as in Figure 11.
\begin{figure}[h]
 \centering
  \includegraphics[width=.45\linewidth]{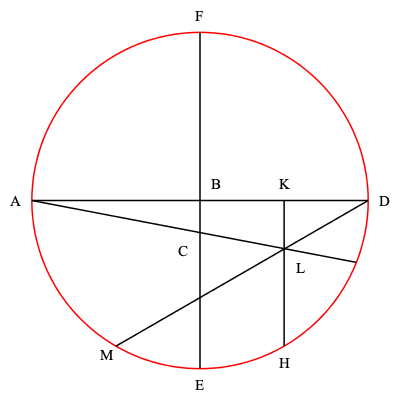}
  \caption{Configuration for finding two mean proportionals between $AB$ and $BC$.}
\end{figure}
Provided that we can find a point $H$ and the vertical segment $KH$ (parallel to $EF$) so that the intersection 
$L$ of the extended line $AC$ and $KH$ makes
the arcs $EH$ and $EM$ (formed by the line through $D$ and $L$) equal, then the ``Lemma of Diocles'' will imply that 
$KH$ and $DK$ are two mean proportionals between $AK$ and $KL$.  But the triangles $ABC$ and $AKL$
are similar, and hence we can rescale all four lines by the ratio $AB : AK$ to get two mean proportionals between 
$AB$ and $BC$ as desired.  

Finding the required point $H$ could be done by the same sort of approximate trial-and-error processes we saw in the previous
methods.  But Diocles and Clavius now actually take this idea one step farther.  Namely, start by considering the 
circle with radius $AB$ as before.   If the locus of all points $L$ as in this 
figure for all possible arcs $EM$ is considered, the so-called \emph{cissoid} of Diocles (a cuspidal cubic algebraic
curve) is obtained.\footnote{Clavius does not use this name, though.  In the coordinate system suggested by placing
the diameters along the coordinate axes and
taking the circle to have radius $1$, the equation of the cissoid is $(x^2 + (y+1))^2 (y+1) = 2x^2$.}  
Namely, for each possible $M$ in the quadrant $AE$, consider the line $DM$ and then take $K$ so that the vertical line $KH$ makes
the arcs $EH$ and $EM$ equal.  Take the point $L$ corresponding to that choice of $M$ as the intersection of the lines
$KH$ and $MD$.  

Then for each point $C$ on the radius $BE$, the line $AC$, when extended, will intersect the cissoid at some uniquely 
determined $L$ and hence produce a line $KH$ making the arcs $EM$ and $EH$ equal.  
Then two mean proportionals between $AB$ and $AC$ will be found as above by rescaling $KH$ and $DK$, which 
are mean proportionals between $AK$ and $KL$.
Thus the cissoid in effect solves the problem for all possible pairs of the fixed $AB$ and smaller segments $BC$ 
\emph{simultaneously}.  As usual, Clavius provides a much more specific description of how the cissoid curve, or more 
precisely, as many points on the curve as desired, can be constructed.  Absent the whole curve, that is having only a 
finite collection of points on the curve, some approximation or judgment of the geometer would still be needed to connect 
the points into a continuous curve and find an appropriate point $L$ for an arbitrary given line $BC$ as above.  Hence the 
qualifier ``approximately'' (the \emph{prope verum} in the Latin) still applies.

The final method for constructing two mean proportionals between given line segments addressed by Clavius is
the one attributed to Nicomedes, using the \emph{conchoid} curve.  (The name was apparently suggested by the similarity between
its shape and the shells of some marine molluscs.)  This discussion is probably the closest Clavius comes
to simply reproducing what he finds in Eutocius or parts thereof.  Clavius starts by saying the conchoid can be drawn with a certain instrument 
(which is described in the first section of Eutocius's version of this method).  But since Clavius does not have a copy of the instrument, 
he says it will be enough to give a construction by which as many points on the conchoid as desired can be produced. (Note the 
parallel with the discussion of the cissoid.)  So let $AB$
be a line and let $CD$ be another perpendicular line meeting $AB$ at a right angle at $E$.  Taking $D$ as a \emph{pole}, consider all
straight lines passing through $D$.  All lines except the parallel to $AB$ through $D$ will intersect $AB$ (extended if necessary).  
Say the line $DS$ meets $AB$ at $S$.  Then extending the line again in the direction of $S$, there will be another point $F$
on the line with $SF = EC$.  The locus of all such points $F$ is the curve known as the conchoid.\footnote{More precisely, if we 
introduce coordinates placing the $x$-axis along the line $AB$ and $E$ at the origin and take $CE = ED$, Nicomedes' conchoid 
is one of the connected components of the real algebraic quartic curve defined by $(x^2+(y+1)^2)y^2 = (y+1)^2$.
There is also a second connected component  below the line $AB$ with a cusp at the point $(0,-1)$, namely the point $D$.}
\begin{figure}[h]
 \centering
  \includegraphics[width=.6\linewidth]{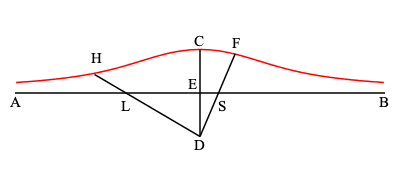}
  \caption{The red curve is the conchoid of Nicomedes.}
\end{figure}
Next, following Eutocius, Clavius proves two ``remarkable properties''\footnote{duas proprietates huius lineae insignes, 
\cite[p.~270]{Clavius1606}.} considered by Nicomedes.  First, the farther the point $S$ is from $E$, the smaller the vertical distance is
from $F$ to the line $AB$ and second, the conchoid meets every line lying above $AB$, no matter how close.\footnote{In modern terms, the line 
$AB$ is a \emph{horizontal asymptote} of the conchoid.}  Following this Clavius shows how the conchoid gives a solution of 
the following problem that Eutocius credits to Nicomedes:
\begin{quote}
Given any rectilinear angle, and a point outside the lines making up 
the angle, to construct from this point, a line intersecting the 
lines containing the given angle, so that the portion of the line
intercepted between the lines is equal to a given line.\footnote{ Dato quouis angulo
rectilineo, \& puncto extra lineas angulum datum comprehendentes: Ab illo
puncto educere rectam secantem rectas datum continentes angulum, ita
vt eius portio inter illas rectas intercepta aequalis sit datae rectae.} (\cite[p.~272]{Clavius1606})
\end{quote}
Finally, Eutocius and Clavius show how the solution of this problem lets one construct the triangle $GDF$ in Figure 9 for 
which $AF$ and $CG$ are the two mean proportionals between the sides $AB$ and $BC$ of the rectangle as in that
figure.  There are several additional constructions of lines made starting from the rectangle and the problem above is used to 
produce a line intersecting two other lines such that the line intercepted is equal to one half of $AB$.   Here Clavius
adds a sort of mnemonic diagram intended to help the reader visualize some of the proportionalities between sides
of similar triangles in the rather complicated figure.  

These discussions give additional excellent examples of the dual practical and theoretical focus of Clavius's \emph{Geometria
Practica} that we identified in the Introduction.  They show how Clavius engages with sources from 
ancient Greek geometry and how he seeks to adapt the results for practical purposes, while still providing a
complete development of the theory involved.  This complete development often includes added or modified 
features designed to smoothe the way for a student learning the material.  In particular, the choices he makes of 
which methods to include certainly do address his criteria of \emph{usefulness}, \emph{ease of application}, and 
\emph{lower susceptibility to error}.  Moreover, the methods of Diocles and Nicomedes are certainly more involved than the previous ones, 
so there is a very clear progression from simpler methods to more complicated ones.

\section{Clavius's presentation of extraction of $n$th roots in Book VI.}

While Clavius was a strong adherent and proponent of the geometrical methods found 
in Euclid's \emph{Elements}, he understood very well that many geometrical constructions corresponded
to \emph{algebraic or numerical} calculations.  A key example is the construction of 
one or two mean proportionals between two lines $AB$ and $CD$.  If $XY$ is one mean proportional, then 
$AB : XY :: XY : CD$, and in algebraic terms 
$$\frac{AB}{XY} = \frac{XY}{CD}.$$
Hence 
$$\left(\frac{XY}{CD}\right)^2 = \frac{AB}{CD},$$
so, in numerical terms, finding $XY$ is essentially the same thing as finding a \emph{square root}.
We have already seen that finding two mean proportionals is essentially the same thing as finding a \emph{cube root}.
The pattern would continue: if any number $n \ge 1$ of mean proportionals in continued proportion were found, that
would be the same as finding an $(n+1)$st root.  As a result many texts on practical geometry, including the texts
of Fibonacci and Pacioli mentioed earlier, included extensive
discussions of numerical algorithms for computing square and cube roots (at least).  
Clavius's book is no exception.  He points out this connection in Proposition 18 immediately following the material
discussed in \S 4 above and he devotes the final section of his Book VI to this topic, starting with the statement 
of Proposition 19:  ``To extract a root of any sort.''\footnote{Radicem cuiuslibet generis extrahere.}

As we mentioned in the Introduction, this is another case where Clavius does not acknowledge a source explicitly.
Indeed, he is almost coy about this, saying only that his treatment of a root extraction algorithm is ``from a book of a certain 
remarkable German arithmetic.''\footnote{ex libro eximij cuiusdam Arithmetici Germani}   The first treatment of this 
material in the German language was contained in the very well-known book \emph{Die Coss} by Christoph Rudolff (1499--1545).\footnote{The 
word \emph{Coss} in German was borrowed from the Italian \emph{cosa} (i.e. ``thing'').  Both were used to represent the unknown in an
algebraic equation before the development of symbolic forms of algebraic expressions.  A generation of 
early German algebraists were known as \emph{cossists}.}
A revised and much expanded edition of this book prepared by Michael Stifel (1487--1567) was published in 
1553\footnote{The title page says, in part: Die Coss Christophs Rudolffs Mit sch\"onen Exempeln der Coss Durch Michael
Stifel Gebessert und sehr gemehrt ... Zu K\"onigsperg in Preussen, Gedr\"uckt durch Alexandrum Lutomyslensem im jar 1553.
(That is, ``Christoph Rudolff's \emph{Die Coss}, with beautiful examples [of these techniques], improved and much augmented by 
Michael Stifel, ... .)'' Digitized version from {\tt www.math.uni-bielefeld.de/\~{}sieben/rudolff.pdf}.} and went through several later editions.
It seems very probable that this (or possibly a later edition) is the book Clavius was drawing from, and the specific section he was
looking at was the  \emph{Anhang} to Chapter 4 of Part I written by Stifel, found starting on folio 46 and going to folio 59.

As usual, Clavius's account is \emph{not directly copied} from Stifel's version.  Clavius's explanations are rewritten and expanded.
Different numerical examples are presented.  Our conjecture that Clavius was consulting this source is based on the fact 
that the overall outline of the method Clavius presents is essentially exactly the same as what Stifel presents:
\begin{itemize}
\item Very similar terminology for the different species of roots is used, e.g. ``zensizenic'' roots are fourth roots, 
``surdesolidic'' roots are fifth roots, and so forth.  Variations of this terminology are found in many 16th century 
works dealing with algebra, though, so this is only a start.
\item The digits from the number whose root is being found are grouped into ``points'' in the same way by marking 
certain digits with dots; each ``point''  will
yield one digit of the root (Clavius writes
the dots below the corresponding digits, while Stifel writes them above, though).
\item Tables of $n$th powers of the digits $1,2,3,\ldots,9$ are provided for use with each ``point'' so that the 
largest $n$th power that can be subtracted from the ``point'' can be identified.
\item The essential role of collection of ``special numbers'' for each species of root 
to be used in preparing the ``divisor'' at each step of the algorithm is the same in both.\footnote{ In our terms, these ``special numbers''
are binomial coefficients times powers of $10$, since the algorithm works with numbers of the form $(10 d_k + d_{k+1})^n$,  
where $d_k$ and $d_{k+1}$ are successive digits of the root.  For
example, 
$$(10 d_k + d_{k+1})^3 = 1000 d_k^3 + 300 d_k^2 d_{k+1} + 30 d_k d_{k+1}^2 + d_{k+1}^3,$$
so the coefficients $300$ and $30$ are the ``special numbers'' used to compute cube roots.} 
\item The calculations are laid out in a very similar (and also very easy) tabular format.
\end{itemize}
 Clavius provides a table containing the binomial coefficients up to $n = 17$ on page 278 of \cite{Clavius1606}.  
 This is not found explicitly in Stifel's discussion, 
 so Clavius might be taking this from another source he does not mention, or computing the entries himself.\footnote{If 
 this the table to which Knobloch is referring at \cite[p.~272]{Knobloch1995}, then the source may be Cardano.}

Probably the best way to convince the reader of this identification of Clavius's source is to quote from two extractions of cube roots, 
one from Clavius and one from Stifel's \emph{Anhang}.  The process described finds the root decimal
digit by decimal digit.  The steps all follow the same pattern after the determination of the left-most digit.
So the point will be made if we look at the determination of the first two digits of the root in the examples.
We begin with the first two steps of this example from Clavius:\footnote{Sit ex numero 
\vskip 5pt
\begin{center}
\begin{tabular}{ccccccccc}
2 &3& 9& 4& 8& 3& 1& 9& 0  \\
   &&$\bullet$&&&$\bullet$&&&$\bullet$
\end{tabular}
\end{center}
\vskip 5pt
extrahenda radix cubica.

Primvm
ex puncto $239.$  subtraho cubum $216.$  qui est maximus in eo 
\[
\begin{array}{rcr}
36&--&300\\
  6&--&30
\end{array}
\]
contentus, cuius radicem $6.$  scribo in Quotiente ad marginem.  Et quia relinquitur numerus                          .
$23.$  erit sequens punctum $23483.$             
Deinde paro diuisorem hoc modo.  Supra radicem inuentam $6.$  pono eius quadratum
$36.$   Et ad dextram colloco duos numeros peculiares radicis cubicae, nimirum $300.$  \&
$30.$  vt hic vides.  Multiplico superiores duos numeros $36.$  \& $300.$  inter se, \& producto
$10800.$  addo productum $180.$  ex multiplicatione numerorum inferiorum $6.$  \& $30.$  
inter se.  Nam summa $10980.$  erit Diuisor.  Satis etiam esset productus ex duobus
superioribus inter se multiplicatis, nimirum $10800.$  pro Diuisore.  quod in alijs extractionibus
intelligendum quoque est.  Diuido ergo punctum meum $23483.$  per diuisorem inuentum 
$10980.$  \& Quotientum $2.$  scribo post figuram $6.$  prius inuentam.  

Pingo post haec figuram huiusmodi.  Ad dextram numerorum $36.$  \& $300.$ 
\[
\begin{array}{rcrcr}
36&--&300&--&2.\\
6&--&30&--&4.\\
&&&& 8.
\end{array}
\]
 colloco inuentam figuram $2.$  \& infra eam eius quadratum $4.$  \& sub
hoc cubum eiusdem $8.$   Nam si tam superiores tres numeri $36.$  $300.$  \& $2.$  quam inferiores
tres $6.$  $30.$  \& $4.$  inter se multiplicentur, \& productis $21600.$  \& $720.$   addatur cubus $8.$  fiet
numerus $22328.$  quem si ex meo puncto $23483.$  subtraham, remanent $1155.$  atque adeo
puncto sequens erit $1155190.$}

\begin{quote}
Let it be required to extract the cube root of 
\vskip 5pt
\begin{center}
\begin{tabular}{ccccccccc}
2 &3& 9& 4& 8& 3& 1& 9& 0  \\
   &&$\bullet$&&&$\bullet$&&&$\bullet$
\end{tabular}
\end{center}
\vskip 5pt
\noindent
First, from the point $239$, I subtract the $216$ which is the largest cube contained in it.
I write its cube root $6$ in the margin in the quotient.  And since $23$ is left over, the next
point will be $23483$.
\[
\begin{array}{rcr}
36&--&300\\
  6&--&30
\end{array}
\]
Next I provide a divisor in this way.  Over [the digit] 6 of the root found above, I put its
square, $36$. And on the right, I place the two particular numbers for cube roots, namely,
$300$ and $30$.  I multiply the numbers on the first row, yielding a product of $10800$ and I
add the product from multiplying the two numbers on the second row, $180$.  The sum
$10980$ will be the divisor.  (It would be enough to take the product of the two numbers
on the first row as the divisor, namely $10800$, as must be understood in other root
extractions.)  I divide the point $23483$ by $10980$ and write the quotient $2$ next  to the digit
$6$ found first.  I treat what comes after this digit as follows.  At the right
\[
\begin{array}{rcrcr}
36&--&300&--&2.\\
6&--&30&--&4.\\
&&&& 8.
\end{array}
\]
of the numbers $36$ and $300$, I add this digit $2$ [found in the quotient] and below it, its square, $4$, and its
cube, $8$.  Now, the three numbers on each of the first two rows are multiplied, 
and the products are $21600$ and $720$.  Adding the cube $8$ makes $22328$.  I subtract
this from the point, leaving $1155$, and the next point will be $1155190$.   (\cite[pp.~280--281]{Clavius1606})
\end{quote}
Clavius continues to find the (approximate) cube root $621$ for $239483190$.  Note that 
$621^3 = 239483061$, so this value is $129$ ``short.''  Later in this section, Clavius
also shows how to compute additional decimal digits in the fractional part, obtaining 
closer approximate cube roots.  

We now translate a step of the computation from folios 47-49 in Stifel's \emph{Anhang}:\footnote{
Exemplum.  
\vskip 5pt
\begin{center}
\begin{tabular}{cccccccccccccccccccc}
 &$\bullet$&&&$\bullet$&&&$\bullet$&&&$\bullet$\\
8 & 0 & 6 & 2 & 1 & 5 & 6 & 8 & 0 & 0 & 0
 \end{tabular}
\end{center}
\vskip 5pt
\noindent
Erstlich subtrahir ich von dem hindersten puncten (das ist von $80$) die 
aller gr\"oste cubic zal/ die ich subtrahiren kan.  Die selbig ist $64.$ so 
bleybett nach vbrig davon $16$ die geh\"oren denn sum nehisten 
puncten hernach/ der selbig uverkompt denn dise figuren $16621$.  So 
setz nu die cubic w\"urzel von $6$ in den quotient.  facit $4.$ und is also 
der erst punct aufsgericht.  

So nehme ich nu fur mich den andern punct/ nemlich $16221.$  Den dividir
ich mit $4800.$ (das kompt von $300$ mal $16$)  Nu gibt das gedacht 
dividiren nur $3$ in den quotient.  Und ist also die newe figur gefunden.

Dem selbigen nach stehn die zwo zalen $300$ und $30$. mit jren zugethonen
zalen also.
\[
\begin{array}{rcrcr}
16&--&300&--&3\\
4&--&30&--& 9\\
\end{array}
\]
Denn erstlich ist gefunden in den quotient de figur 4. die steht neben $30$ 
zur lincken hand/ vnd drob neben $300$ steht jr quadrat/ nemlich $16$.

So is nu darnach gfunden in den quotient die figur $3$. Die steht oben 
neben $300$ zur rechten hand/ vnd darunder steht jr quadrat $9.$ neben
$30.$  wie du alles vol sihest.  

So multiplicir ich nu/ vnd sprich. $16$ mal $300$ mal $3$. facit. $14400$.
vnd $4$ mal $30$ mal $9.$ facit $1080$.   Das addir ich/ so kompt 
$15480$.  Das subtrahir ich von $16621$.  Als vom andern puncten
diser operation/ so bleyben denn $1141$.

Auffs letzt multiplicir ich die newe gefundne figur Cubice.  Nemlich
$3$ mal $3$ mal $3.$ facit $27$.  die subtrahir ich auch/ so bleyben
$1141$. die geh\"oren zu volgenden punct.}

\begin{quote}
Example.  
\vskip 5pt
\begin{center}
\begin{tabular}{cccccccccccccccccccc}
 &$\bullet$&&&$\bullet$&&&$\bullet$&&&$\bullet$\\
8 & 0 & 6 & 2 & 1 & 5 & 6 & 8 & 0 & 0 & 0
 \end{tabular}
\end{center}
\vskip 5pt
\noindent
First I subtract the largest cube that I can from the leftmost point (that is, from $80$).
That is $64$, leaving $16$, which then belongs to the next point, which is
composed of the digits $16621$.  So now I set the cube root of $64$ in the quotient, 
and the first point is decided.  So then I take the next point, namely $16621$.  
I divide that by $4800$ (which is $300$ times $16$) and that division gives $3$
in the quotient.  And so the new digit is found.  I put this next to the two numbers
$300$ and $30$ with the accompanying numbers in this way:
\[
\begin{array}{rcrcr}
16&--&300&--&3\\
4&--&30&--& 9\\
\end{array}
\]
Since the digit 4 was found [first] in the quotient, that is placed next to the $30$ on 
the left, and above, next to $300$ goes its square, namely $16$.  On the right 
next to the $300$ goes the next digit $3$, and its square $9$ goes
below next to the $30$, as you clearly see.

So now I multiply and say $16$ times $300$ times $3$ makes $14400$, and
$4$ times $30$ times $9$ makes $1080$.  I add those and obtain $15480$.
I subtract that from the $16621$ as from the other points.  The number $1141$
remains.  Last, I take the cube of the newly-found digit $3$.  Namely,
$3$ times $3$ times $3$ makes $27$. I also subtract this and $1114$ remains.
This belongs to the following point.
\end{quote}
Since there were four ``points'' in the original number, Stifel's cube
root will contain four decimal digits.  After two more steps of the process,
he finds the value $4320$, an exact cube root of 
$80621568000$.

If my conjecture that Clavius was following Stifel's presentation of a root extraction 
algorithm here is correct (and I hope I have proved the point with the quotations above!), 
then there remains the question why Clavius did not
make an explicit attribution to Rudolff and/or Stifel.
It is certainly possible that Clavius thought he did not need to say any more to identify the 
source because Rudolff's \emph{Die Coss} was extremely well-known, at least in 
German-speaking areas because it was the 
first book on this material published in German.  However, there is another circumstance that
might just provide another component of an explanation.  Namely, Stifel had
started out as an Augustinian monk, but later became a Protestant minister and an outspoken supporter
of Martin Luther.  Unlike the citation of the Protestant Scaliger mentioned earlier in \S 4, where 
Clavius was being explicitly critical of the other's work,  Clavius was singling out
Stifel's algorithm for high praise and recommending its use. 
Under those circumstances, it may be that Clavius (or his Jesuit colleagues and superiors) thought it was not politic to 
mention Stifel's name.    

\section{Conclusions}

Clavius presented a tremendous amount of interesting and useful mathematics in his
\emph{Geometria Practica} and in his other writings.  In assembling the material for 
this book, he drew on an extremely broad range of ancient, medieval, and contemporary 
sources.  At the same time, his typical procedure was to rework, augment, and clarify the 
mathematical texts he dealt with.  It seems arguable that he achieved his stated goal of presenting
the whole range of practical geometry as understood in his time, and he did it in a form
that would be useful for his readers.  

The quality of this work was recognized very soon after it
appeared, as evidenced (for instance) by the fact that mathematicians such as Kepler
mentioned sections of this book in their writings.\footnote{In Kepler's \emph{Harmonice Mundi}.
The reference was to Clavius's discussion of various approximate constructions of regular heptagons 
in Book VIII.}  Recognition of Clavius's work was also evident in other ways.  In 
the Jesuit mission in China, one of Clavius's former students in the \emph{Collegio Romano}, Mateo Ricci, 
S.J. (1552--1610), together with his Chinese collaborator Xu Guangqi (1562--1633),
made translations of not only the first six books of Euclid's \emph{Elements} from 
Clavius's version, but also material from the \emph{Geometria Practica}.  
Later, Giacomo Rho, S.J. (1593--1638) made Chinese
translations of additional sections of this work.\footnote{\cite[pp.~318--319]{Martzloff})}

However, if I may be allowed to speculate in this last paragraph, in some ways, 
I would argue that Clavius's \emph{Geometria Practica} actually represents almost the end of 
the sub-genre of \emph{``theoretical practical geometry''} in the style we have seen in our discussions.
There were certainly many later practical geometry books, but they tended more toward the ``practical,'' and
in many cases omitted proofs or theoretical developments.  In addition, Clavius's essential mathematical
conservatism and his devotion to the synthetic Euclidean tradition in geometry would shortly come to seem 
very old-fashioned. The recovery of Pappus's treatment of the Greek tradition of geometric \emph{analysis} in 
Book VII of the \emph{Mathematical Collection},
combined with the ever-growing influence of algebraic thinking was the impetus for an explosion of work 
starting in the late 16th century and continuing into the first half of the 17th 
century (this is discussed in a fascinating way in \cite{Bos}).  
 But this was largely orthogonal to the ways that Clavius approached geometry and he seemingly had
 little interest in or taste for that side of Pappus's writings.\footnote{See in particular \cite[Chapter 2, \S 3]{Sasaki}.}
Within 30 years of his death, the introduction and systematic use of \emph{analytic, or coordinate geometry} 
by Descartes and others was well under way.  That new way of harnessing the power of algebra to \emph{discover new
geometrical results} and prove them was fundamentally changing the practice of mathematics. 


\bibliographystyle{plain}

\end{document}